\documentclass{article}
\usepackage{amsmath, amssymb}
\usepackage{amssymb,pb-diagram,pst-all,graphicx}
\usepackage{amscd}

\usepackage{fancybox}

\newcommand{\ZZ}{\mathbb Z}
\newcommand{\PP}{\mathbb P}

\newcommand{\CC}{\mathbb C}

\newcommand{\mcB}{\mathcal B}
\newcommand{\mcC}{\mathcal C}

\newcommand{\mcL}{\mathcal L}

\newcommand{\mcD}{\mathcal D}
\newcommand{\mcQ}{\mathcal Q}
\newcommand{\mcT}{\mathcal T}

\newcommand{\bbA}{\mathbb A}

\newcommand{\tor}{{\mathrm {tor}}}

\newcommand{\Cov}{\mathop {\rm Cov}\nolimits}

\newcommand{\MW}{\mathop {\rm MW}\nolimits}

\newcommand{\NS}{\mathop {\rm NS}\nolimits}

\newcommand{\Red}{\mathop {\rm Red}\nolimits}

\newcommand{\Sing}{\mathop {\rm Sing}\nolimits}

\def\miya#1{\mathrel{\mathop{\rightarrow}\limits^#1}}

\newcommand{\mcE}{\mathop {\cal E}\nolimits}

\newcommand{\mcO}{\mathop {\mathcal O}\nolimits}

\newtheorem{thm}{Theorem}[section]

\newtheorem{prop}{Proposition}[section]
\newtheorem{lem}{Lemma}[section]
\newtheorem{defin}{Definition}[section]
\newtheorem{exmple}{Example}[section]
\newtheorem{rem}{Remark}[section]
\newtheorem{qz}{Question}[section]

\newcommand{\I}{\mathop {\rm I}\nolimits}
\newcommand{\II}{\mathop {\rm II}\nolimits}
\newcommand{\III}{\mathop {\rm III}\nolimits}

\newenvironment{remark}{\begin{rem}\rm }{\end{rem}}

\newcommand{\qed}{\hfill $\Box$}

\newcommand{\proof}{\noindent{\textsl {Proof}.}\hskip 3pt}

\renewcommand{\thesubparagraph}{\theparagraph.\@arabic\c@subparagraph}
\begin{document}
  
  \begin{center}
{\Large \bf Rational points of elliptic surfaces 
and 
Zariski $N$-ples for cubic-line, cubic-conic-line arrangements
}

\bigskip

\bigskip
{\bf 
Shinzo BANNAI, Hiro-o TOKUNAGA and  Momoko YAMAMOTO
}

\end{center}
\normalsize

\begin{abstract}
In this paper, we continue the study of the relation between rational points of rational elliptic surfaces and plane curves. As an application, we give first examples of Zariski pairs of cubic-line arrangements that do not involve inflectional tangent lines.
\end{abstract}

\section*{Introduction}

  In this article, we study the arithmetic of  rational points of certain rational elliptic surfaces from a geometric point of view 
  in order to construct arrangements of plane curves of low degree which give  rise to candidates for Zariski pairs.
  A pair of reduced plane curves
  $(\mcB^1, \mcB^2)$ is said to be a {\it Zariski pair}
 if it satisfies the following conditions:

\begin{enumerate}

\item[(i)] For each $i$, there exists a tubular neighborhood $T(\mcB^i)$ of $\mcB^i$ such
that $(T(\mcB^1), \mcB^1)$ is homeomorphic to $(T(\mcB^2), \mcB^2)$.

\item[(ii)] There exists no homeomorphism from 
$(\PP^2, \mcB^1)$ to $(\PP^2, \mcB^2)$.
\end{enumerate}

The first condition can be replaced by {\it the combinatorics} (or {the combinatorial type}) of $\mcB^i$. For the precise definition
of the combinatorics,  see \cite{act} (It can also be found in \cite{tokunaga14}).  Since the combinatorics is easier to
treat with, we always consider that of $\mcB^i$.
The study of Zariski pairs was originated by Zariski in \cite{zariski29}. Since the '90's there have been
a lot of results on Zariski pairs by many mathematicians via various methods (see the reference \cite{act},
for example).
As we noted in \cite{act}, there are two main ingredients in the study of Zariski pairs.
Namely,
\begin{enumerate}

\item[(I)] To find reduced curves $\mcB^1$, $\mcB^2$ with the same combinatorics so that $\mcB^1$, $\mcB^2$ have
certain different features, and

\item[(II)] To prove $(\PP^2, \mcB^1)$ is not homeomorphic to $(\PP^2, \mcB^2)$ based on the different feature as above.

\end{enumerate}

A Zariski $N$-ple is a natural generalization, where the number of curves is increased.
One of our new feature of this article concerns   (I): Construction of plane curves via geometry and arithmetic 
of sections for certain rational elliptic surfaces. This basic idea can be found in
\cite{tokunaga14} by the first autor  and in \cite{bannai-tokunaga}  by the first and second authors. 
In this article, however,
we make use of the arithmetic of sections more intensively than previous papers.

For (II), in order to distinguish $(\PP^2, \mcB^1)$ and $(\PP^2, \mcB^2)$, our tool  is Galois covers
 branched along $\mcB^i$  developed in \cite{bannai-tokunaga, tokunaga14}.  
Note that there are
various other tools, for example, 
the fundamental group $\pi_1(\PP^2\setminus\mcB^i, \ast)$, braid monodromy and Alexander invariants.
Recently two more new tools, the linking set and the connected number,  are introduced
by J.-B. Meilhan, B.~Guerville-Ball\'{e} \cite{benoit-jb} and T.~Shirane \cite{shirane16}, respectively.

 We explain our object more concretely.
  The first and second authors have studied Zariski pairs (or $N$-ples) for arrangements of curves whose
 irreducible components are of low degree, less than or equal to 4, via geometry of sections and bisections of
 rational elliptic surfaces (\cite{bannai-tokunaga, tokunaga14}).  In this article, we also continue to study such objects. More precisely, we consider
 a reducible curve whose irreducible components consist of 
 
 \begin{enumerate}
 \item[(i)] one irreducible cubic and lines, and
 
 \item[(ii)] one irreducible cubic, smooth conics and lines.
 \end{enumerate}

In \cite{artal94}, a Zariski pair of sextics consisting of 
 a smooth cubic and its three inflectional tangents was given. This example
was also considered in \cite{tokunaga96} from a different approach. In \cite{benoit-jb}, a Zariski
pair of octics consisting of a smooth cubic and its $5$ inflectional tangents is given. In \cite{bbst},  
that of reducible curves consisting of 
a smooth cubic and  $k$ of its inflectional tangents ($k = 4, 5, 6$) are considered. Note that
there exist no Zariski pair of a smooth cubic and $k$ of its inflectional tangents for $k = 1, 2, 7, 8, 9$.
Also, in \cite{at00}, E.~Artal Bartolo and the second author studied a Zariski pair of sextics whose
irreducible components are a nodal cubic, a smooth conic and an inflectional tangent line.
All of these examples  contain inflectional tangents of a cubic as irreducible components.
As for another new feature of this article,  we focus on  Zariski pairs for cubic-line or cubic-conic-line arrangements 
{\it without inflectional tangents}. Also, since Zariski pairs for sextic curves have already
been intensively studied by many mathematicians, e.g., \cite{oka05, oka06, shimada}, we consider
the case of degree $7$ as follows:


{\bf Combinatorics 1.} Let $\mcE$, $\mcL_o$ and $\mcL_i$ ($i = 1, 2, 3$) be as below and
we put $\mcB = \mcE  + \mcL_o + \sum_{i=1}^3 \mcL_i$:

\begin{enumerate}

\item[(i)] $\mcE$:  a smooth or nodal cubic curve.

\item[(ii)] $\mcL_o$: a transversal line to $\mcE$ and we put $\mcE \cap \mcL_o = \{p_1, p_2, p_3\}$.

\item[(iii)] $\mcL_i$ ($i = 1, 2, 3$): a line through $p_i$ and tangent to  $\mcE$ at a different point from $p_i$.
We denote the tangent point of $\mcL_i$ by $r_i$.

\item[(iv)] $\mcL_1, \mcL_2$ and $\mcL_3$ are not concurrent.

\end{enumerate}


{\bf Combinatorics 2.} Let  $\mcE$, $\mcC$, $\mcL_o$ and $\mcL$ as below and
we put $\mcB = \mcE + \mcC +  \mcL_o + \mcL$\,:

\begin{enumerate}

\item[(i)] $\mcE$:  a nodal cubic curve.

\item[(ii)] $\mcL_o$: a transversal line to $\mcE$ and we put $\mcE \cap  \mcL_o = \{p_1, p_2, p_3\}$.

\item[(iii)] $\mcL_i$: a  line connecting the node of  $\mcE$  and one of $p_i, (i = 1, 2, 3)$.

\item[(iv)] $\mcC$: a contact conic  to $\mcE + \mcL_o$ and intersecting $\mcL$ transversely.

\end{enumerate}

 Here
 we call a smooth conic $\mcC$ a contact conic to a reduced plane curve $\mcB$ if the following condition is satisfied:
$(\ast)$ Let $I_x(\mcC, \mcB)$ denotes the intersection multiplicity at $x \in \mcC \cap \mcB$.
For $\forall x \in \mcC \cap \mcB$, $I_x(\mcC, \mcB)$ is even and $\mcB$ is smooth at $x$. 

In both combinatorics, no inflectional tangents are involved and this is the new feature compared to 
previous examples. 
Let us explain more precisely. In the following, we use
the  notation and terminology  used in \cite{bannai-tokunaga}.

%
%
%
%

Put $\mcQ = \mcL_o + \mcE$ and choose a smooth point $z_o \in \mcE$. 
Consider the minimal resolution $S_{\mcQ}$ of a double cover of $\PP^2$ branched along $\mcQ$ and blow up $S_{\mcQ}$  twice at the 
preimage of $z_o$. Then we obtain a rational elliptic surface 
$\varphi_{\mcQ, z_o} : S_{\mcQ, z_o} \to \PP^1$ and its generic fiber is denoted by $E_{\mcQ, z_o}$(see \S 1 for a more precise description). 
 We denote the induced generically $2$-to-$1$ morphism from $S_{\mcQ, z_o}$ to $\PP^2$ by $f_{\mcQ, z_o}$. Let  $E_{\mcQ,z_o}$ be
 the generic fiber of $\varphi_{\mcQ, z_o} : S_{\mcQ, z_o} \to \PP^1$.  The group of sections of $\varphi_{\mcQ, z_o}$, 
 $\MW(S_{\mcQ, z_o})$,  can be canonically  identified with the group of $\CC(t)$-rational points $E_{\mcQ, z_o}(\CC(t))$.
  For $s \in \MW(S_{\mcQ, z_o})$, we denote the corresponding rational  point by $P_s$. Conversely, for $P \in E_{\mcQ, z_o}(\CC(t))$,
 we denote the corresponding section by $s_P$. Now, since $s \in \MW(S_{\mcQ, z_o})$ can be considered a curve,
$f_{\mcQ, z_o}(s)$ gives rise to a plane curve in $\PP^2$. In our construction of plane curves with Combinatorics 1 and 2,
we make use of lines and conics of the form $f_{\mcQ, z_o}(s)$ for some $s \in \MW(S_{\mcQ, z_o})$. In our particular cases,  it can be
explained more explicitly as follows (We use the notation in 
\cite{oguiso-shioda} in order to describe the structure of $E_{\mcQ, z_o}(\CC(t))$):
%
%
%
\medskip

{\bf Combinatorics 1.}
(a) $\mcE$: a smooth cubic.  If we choose $z_o \in \mcE \setminus \{p_1, p_2, p_3\}$, 
$E_{\mcQ, z_o}(\CC(t))$ $ \cong D_4^* \oplus \ZZ/2\ZZ$.  As we see in \S 3 , 
we choose generators  $P_0, P_1, P_2, P_3$ for the $D_4^*$-part suitably and the $2$-torsion $P_{\tau}$.    We also 
put $P_4 := P_2 \dot{-}P_3 \dot{+}P_{\tau}$, where $\dot{+}, \dot{-}$ denote the addition and subtraction on $E_{\mcQ, z_o}(\CC(t))$.
 Let $s_{P_i}$ be the sections corresponding to $P_i$ ($i = 1, \ldots, 4$),
respectively.
Then $\mcL_i := f_{\mcQ, z_o}(s_{P_i})$, $(i = 1, \ldots, 4)$ are
tangent lines to $\mcE$ such that each of them passes through $\mcE \cap \mcL_o$.  If $p_i$  $(i = 1, 2, 3)$ are not 
inflection points and any three of them are not concurrent, then
both $\mcB^1:= \mcQ + \sum_{i=1}^3 \mcL_i$ and 
$\mcB^2:= \mcQ + \sum_{i=2}^4 \mcL_i$ have Combinatorics 1-(a).


(b) $\mcE$: a nodal cubic.   If we choose $z_o \in \mcE \setminus \{p_1, p_2, p_3\}$,
$E_{\mcQ, z_o}(\CC(t)) \cong (A_1^*)^{\oplus 3} \oplus \ZZ/2\ZZ$. As we see in \S 3, 
we choose generators  $P_1, P_2, P_3$  of $E_{\mcQ, z_o}$ so that the Gram matrix of
$P_i$ $(i = 1, 2, 3)$ is $[ \langle P_i, P_j \rangle ] = 1/2\delta_{ij}$  and the $2$-torsion $P_{\tau}$.  We put
\[
\begin{array}{ccc}
P_4  :=   P_3 \dot{+} P_1 \dot{+} P_\tau & &P_5  :=  P_1 \dot{+} P_2 \dot{+} P_\tau \\
P_6  :=  P_2 \dot{+} P_3 \dot{+} P_\tau  &  &P_7  :=  P_3 \dot{-} P_1 \dot{+} P_\tau .
\end{array}
\]
Then we infer that $\mcL_i:= f_{\mcQ, z_o}(s_{P_i})\, (i = 4, 5, 6, 7)$ are tangent lines to $\mcE$ such that each of them passes
 through $\mcE \cap \mcL_o$.  If any three of them are not concurrent, then
both $\mcB^1:= \mcQ + \sum_{i=4, 5, 6} \mcL_i$ and 
$\mcB^2:= \mcQ + \sum_{i= 5, 6, 7} \mcL_i$ have Combinatorics 1-(b).

{\bf Combinatorics 2.} We keep our notation in  Combinatorics 1.
Our construction is similar to that in \cite{tokunaga14}.  We first note that $\mcL_i := f_{\mcQ, z_o}(s_{P_i})$ $(i = 1, 2, 3)$ are lines
connecting the node of $\mcE$ and $p_i$.  Plane curves given by  $f_{\mcQ, z_o}(s_{[2]P_i})$  $(i = 1, 2, 3)$ are
contact conics by a similar argument to that in \cite[p. 633]{tokunaga14}.
Suppose that $\mcC_i$ meets $\mcL_j$  transversely for any $i, j$. 
Then $\mcB^1 := \mcQ + \mcL_i +  \mcC_i $ and $\mcB^2:= \mcQ + \mcL_i +  \mcC_j$ ($i \neq j$) have
Combinatorics 2.


\medskip

Now our statement is

\begin{thm}\label{thm:comb-1}{Let $\mcB^1$ and $\mcB^2$ be as above. Then $(\mcB^1, \mcB^2)$ is a Zariski pair.}
\end{thm}

The organization of this paper as follows:
In \S 1, we give a summary on  facts and previous results, which we need later.  
We prove Theorem~\ref{thm:comb-1} in \S 2. A Zariski triple and 4-ple for
cubic-conic-line or cubic-conic arrangements are considered in \S 3 and explicit examples are given in \S 4.

\medskip

{\bf Acknowledgements.} The second author is partially supported by Grant-in-Aid for Scientific Research C (17K05205).


\section{Preliminaries}

\subsection{Elliptic surfaces}

\subsubsection{Generalities}

As for basic references about elliptic surfaces, we refer  to \cite{kodaira, miranda-basic}.  We also
refers to \cite{shioda90} for general facts on the Mordell Weil lattices. In particular, 
for those on rational elliptic surfaces, we refer to \cite{oguiso-shioda}.
We also use the notation and terminology used in \cite{bannai-tokunaga, tokunaga14} freely.
In this article, by an {\it elliptic surface}, we always mean a smooth projective surface $S$ 
  with a fibration $\varphi : S \to C$ over a smooth projective curve $C$ as follows:
 \begin{enumerate}
 
  \item[(i)]  $\varphi$ has a section
 $O : C \to S$ (we identify $O$ with its image). 
  \item[(ii)] There exists a non-empty finite subset, $\Sing(\varphi)$, of $C$ such
  that $\varphi^{-1}(v)$ is a smooth curve of genus $1$ (resp. a singular curve ) for $v \in C\setminus \Sing(\varphi)$ (resp. $v \in \Sing(\varphi))$. Note that there exist no multiple fibers since $\varphi$ has the section $O$.

 \item[(iii)]  $\varphi$ is minimal, i.e., there is no exceptional
 curve of the first kind in any fiber. 
 \end{enumerate}
 
 For $v \in \Sing(\varphi)$, we put $F_v = \varphi^{-1}(v)$. 
 We denote its irreducible decomposition by 
 \[
 F_v = \Theta_{v, 0} + \sum_{i=1}^{m_v-1} a_{v,i}\Theta_{v,i}, 
 \]
 where $m_v$ is the number of irreducible components of $F_v$ and $\Theta_{v,0}$ is the
unique irreducible component with $\Theta_{v,0}O = 1$. We call $\Theta_{v,0}$ the {\it identity
 component}.  The classification  of singular fibers is well-known (\cite{kodaira}). We use the Kodaira Notation to denote the types of singular fibers. 
  We also denote the subset  of $\Sing(\varphi)$  consisting of points giving reducible singular fibers by
 $\Red(\varphi) := \{v \in \Sing(\varphi)\mid \mbox{$F_v$ is reducible}\}$. 
Let $\MW(S)$ be the set of sections of $\varphi : S \to C$.  $\MW(S) \neq \emptyset$ as $O \in \MW(S)$.
By \cite[Theorem 9.1]{kodaira}, $\MW(S)$  is an abelian group with the zero
element $O$. We call
$\MW(S)$  the Mordell-Weil group.   
We also denote the multiplication-by-$m$ map ($m \in \ZZ$) on
$\MW(S)$ by $[m]s$ for $s \in \MW(S)$. 
 On the other hand, the generic fiber $E:= S_{\eta}$ of
$S$ is as a curve of genus $1$ over $\CC(C)$, the rational function field of $C$. The restriction of $O$
to $E$ gives rise to a ${\mathbb C}(C)$-rational point of $E$, and one can regard $E$
as an elliptic curve over ${\mathbb C}(C)$, $O$ being the zero element. $\MW(S)$ can be identified
with the set of ${\mathbb C}(C)$-rational points $E(\CC(C))$ canonically. 
For $s \in \MW(S)$, we denote the corresponding rational point by $P_s$. 
Conversely,
for an element $P \in E(\CC(C))$,  we denote the corresponding  section by $s_P$.
  The abelian group $G_{\Sing(\varphi)}$, and the homomorphsim
$\gamma : \MW(S) \to G_{\Sing(\varphi)}$   are those defined in \cite[p. 83]{tokunaga12}.
For $s \in \MW(S)$, $\gamma(s)$ describes at which irreducible component $s$ meets on 
$F_v$.

Let $\NS(S)$ be the N\'eron-Severi group of $S$ and let $T_{\varphi}$ be the 
subgroup of $\NS(S)$ generated by $O$, a fiber $F$ and  $\Theta_{v,i}$ $(v \in \Red(\varphi)$,
$1 \le i \le m_v-1)$. Then we have the following theorems:

\begin{thm}\label{thm:shioda-basic0}(\cite[Theorem~1.2, 1.3]{shioda90})
 {Under our assumptions,

\begin{enumerate}
\item[(i)] $\NS(S)$ is torsion free, and

\item[(ii)] there is a natural map $\tilde{\psi} : \NS(S) \to \MW(S)$ which induces an isomorphism of 
groups
\[
\psi : \NS(S)/T_{\varphi} \cong \MW(S) (\cong E(\CC(C))).
\]
In particular, $\MW(S)$ is a finitely generated abelian group.
\end{enumerate}
}
\end{thm}

\begin{thm}\label{thm:shioda-basic}{(\cite[Theorem~1.3]{shioda90}) Under our assumptions,
there is a natural map $\tilde{\psi} : \NS(S) \to \MW(S)$ which induces an isomorphism of 
groups
\[
\psi : \NS(S)/T_{\varphi} \cong \MW(S).
\]
In particular, $\MW(S)$ is a finitely generated abelian group.
}
\end{thm}

For a divisor on $S$, we  put $s(D) = \tilde\psi(D)$. 

In \cite{shioda90}, a lattice structure on $E(\CC(C))/E(\CC(C))_{\tor}$ is defined by using
the intersection pairing on $S$ through $\psi$. We use the terminologies, notation and results in \cite{shioda90}, freely.
In particular, $\langle \, , \, \rangle$ denotes the height pairing and $\mbox{\rm Contr}_v$ denotes the contribution term
given in \cite{shioda90} in order to compute $\langle \, , \, \rangle$.

\subsubsection{Double cover construction of an elliptic surface}\label{subsec:double-cover}

We refer to \cite[Lectures III and IV]{miranda-basic} for details.
Let $\varphi: S \to \PP^1$ be an elliptic surface over a smooth projective curve $\PP^1$. 
As we see in \cite{bannai-tokunaga, tokunaga14}, $S$ can be represented as the minimal resolution of a double cover of a Hirzebruch surface $\Sigma_d$
as follows. The inversion of $E$ with respect to
the group law  induces an involution $[-1]_{\varphi}$ on $S$. 
Let $S/\langle [-1]_{\varphi}
\rangle$ be the quotient by $[-1]_{\varphi}$.  It is known that  $S/\langle[-1]_{\varphi}\rangle$ 
is smooth and we can  blow down  $S/\langle [-1]_{\varphi} \rangle$ to its relatively minimal
model $W$. We denote the morphisms involved by
      \begin{itemize}
       \item $f: S \to S/\langle [-1]_{\varphi}\rangle$: the quotient morphism, 
       \item $q: S/\langle [-1]_{\varphi}\rangle \to W$: the blow down, and 
      \item $S \miya{{\mu}} S' \miya{{f'}} W$: the Stein factorization of $q\circ f$. 
      \end{itemize}
      Then (i)  $W$ is   $\Sigma_d$, where $d = 2\chi(\mcO_S)$ and (ii)
 the branch locus $\Delta_{f'}$ of $f'$ is of the form $\Delta_0 + \mcT$, where
 $\Delta_0$ is a section with $\Delta_0^2 = -d$ and $\mcT \sim 3(\Delta_0 + d{\mathfrak f})$, 
 ${\mathfrak f}$ being a fiber of the ruling $\Sigma_d \to {\mathbb P}^1$. Moreover, 
 singularities of $\mcT$ are at most simple singularities (see \cite[Chapter II, \S 8]{bpv} for simple singularities and their notation).  
 

Conversely,  if $\Delta_0$ and $\mcT$ on $\Sigma_d$, $d$: even,  satisfy the above conditions, we
obtain an elliptic surface $\varphi : S \to \PP^1$, as the canonical 
resolution of a double cover $f' : S' \to \Sigma_d$ with $\Delta_{f'} = \Delta_0
+ \mcT$, and the following diagram (see \cite{horikawa} for the canonical resolution):
\[
\begin{CD}
S' @<{\mu}<< S \\
@V{f'}VV                 @VV{f}V \\
\Sigma_d@<<{q}< \widehat{\Sigma}_d.
\end{CD}
\]
Here, $q$ is a composition of blowing-ups so that $\widehat{\Sigma}_d
= S/\langle [-1]_{\varphi}\rangle$. 
Hence any elliptic surface  is obtained in this way. In the case when $S$ is rational, $d = 2$.
In the following, we call the diagram above
{\it the double cover diagram for $S$}.

Moreover, if $S$ is rational and has a reducible singular fiber, $\widehat{\Sigma}_2$ can be blown
down to $\PP^2$, as we remark in \cite[1.2.2]{bannai-tokunaga}.  $\mcT$ is mapped to
a reduced quartic $\mcQ$, which in not concurrent four lines,  and $\Delta_0$ is mapped to a smooth point $z_o$ on
 $\mcQ$.  Conversely, given a reduced quartic $\mcQ$ ($\neq$ concurrent four lines) and a point 
 $z_o \in \mcQ$, we obtain $S$ as above, which we denote by $S_{\mcQ, z_o}$, which is nothing but the surface described in the introduction. The induced
 generically $2$-$1$ morphism from $S_{\mcQ, z_o}$ to $\PP^2$ is $f_{\mcQ, z_o}$. 

 \subsubsection {$S_{\mcQ, z_o}$ for the case when $\mcQ$ is $\mcE + \mcL_o$ in the combinatorics 1 and 2}
 
 Let $z_{o}$ be a smooth point on $\mcQ$. The tangent line $l_{z_{o}}$ to $\mcQ$ at $z_{o}$ gives rise to a singular fiber of ${\varphi}_{\mcQ , z_{o}}$ whose type is determined by how $l_{z_{o}}$ intersects with $\mcQ$ as follows:

\begin{center}
\begin{tabular}{|c|c|l|} \hline
(i) & $\I_{2}$ & $l_{z_{o}}$ meets $\mcQ$ with two other distinct points. \\ \hline
(ii) & $\III$ & $z_{o}$ is an inflection point of $\mcE$. \\ \hline
(iii) & $\I_{0}^{\ast}$ & $l_{z_{o}} = {\mcL}_{o}$. \\ \hline
(iv) & $\I_{4}$ & $l_{z_{0}}$ passes through a  point in  $\mcQ\cap \mcL_o$. \\ \hline
\end{tabular}
\end{center}


 Hence by  \cite[Table 6.2]{miranda-persson},  possible configurations of reducible singular fibers of $S_{\mcQ, z_o}$ are
 as follows:

\noindent{\sl Case 1: $\mcE$ is a smooth cubic}
 

\begin{center}
\begin{tabular}{|c|c|} \hline
 & singular fibers \\ \hline
(i) , (ii)& $\{ a \I_{2}, b\III \}, a+ b = 4, a, b \ge 0, b \neq 3$ \\ \hline
(iii) & $\{ I_0^{\ast} \}$ \\ \hline
(iv) & $\{ \I_{4}, 2\I_{2} \}$ \\ \hline

\end{tabular}
\end{center}

%
 
\noindent  {\sl Case 2: $\mcE$ is a nodal cubic}

\begin{center}
\begin{tabular}{|c|c|} \hline
case & configration of singular fibers \\ \hline
(i), (ii)  & $\{a\I_2, b\III\}, a+ b = 5, 0 \le b \le 2$ \\ \hline
(iii) & $\{\I_0^*, \I_2\}, \{\I_0^*, \III\}$ \\ \hline
(iv) &  $\{\I_4, 3\I_2\}$ \\ \hline
\end{tabular}
\end{center}

For our proof of Theorem~\ref{thm:comb-1}, we choose  $z_o$ satisfying  (i) or (ii).
In these cases, since 
the difference between fibers of type $\mathrm{III}$ and $\mathrm{I}_{2}$ do not affect the structure of $E_{\mcQ, z_o}(\CC(t))$ and since  
$\mcL_o$ gives rise to a $2$-torsion section, by \cite{oguiso-shioda}, we infer that
 $E_{\mcQ, z_o}(\CC(t)) \cong D_4^{\ast} \oplus
\ZZ/2\ZZ$ (resp. 
 $(A_1^*)^{\oplus 3} \oplus\ZZ/2\ZZ$) for Case 1 (resp.  Case 2).



%
%

\subsection{Galois covers}

For the notation and terminology  on Galois covers, we use those in \cite{act, bannai-tokunaga, tokunaga94} freely.

\subsubsection{$D_{2p}$-covers}

We here introduce notation for dihedral covers which we use frequently. For details, see \cite{tokunaga94}.
Let $D_{2p}$ be the dihedral group of order $2p$, where $p$ is an  odd prime. 
 In order to present $D_{2p}$, we use
 the notation
 \[
 D_{2p} = \langle \sigma, \tau \mid \sigma^2 = \tau^p = (\sigma\tau)^2 = 1\rangle.
 \]
 Given a $D_{2p}$-cover, we obtain a double cover $D(X/Y)$ of $Y$ canonically by considering the
 $\CC(X)^{\tau}$-normalization of $Y$, where $\CC(X)^{\tau}$ denotes the fixed field 
 of the subgroup of $D_{2p}$ generated by $\tau$.  Then $X$ is a $p$-fold cyclic cover of $D(X/Y)$ and
 we denote the covering morphisms by
 $\beta_1(\pi) : D(X/Y) \to Y$ and $\beta_2(\pi) : X \to D(X/Y)$, respectively.  

\subsubsection{Elliptic $D_{2p}$-covers}

Let $\varphi : S \to \PP^1$ be a rational elliptic surface and let $f : S \to \widehat{\Sigma}_d$ denote the one 
int the double cover diagram.  For our criterion to distinguish  the topology of plane curves, we make use  of the existence/non-existence
 of $D_{2p}$-covers $\pi_p : X_p \to \widehat{\Sigma}_d$ of $\widehat{\Sigma}_d$ satisfying  (i) $D(X_p/\widehat{\Sigma}_d) = 
 S$ and (ii) $\beta_1(\pi_p) = f$.  Following \cite{tokunaga14}, we call such a $D_{2p}$-cover  an
 elliptic $D_{2p}$-cover. We denote the covering transformation of $f$ by $\sigma_f$. 
 As we remark in \cite[\S 3]{tokunaga14}, the branch locus $\Delta_{\beta_2(\pi_p)}$ of 
 $\beta_2(\pi_p)$ is the form
 \[
 \mcD + \sigma_f^*\mcD + \Xi + \sigma_f^*\Xi,
 \]
 where
 \begin{enumerate}
  
  \item[(i)]  no irreducible component of $\mcD$ and $\sigma_f^*\mcD$ is contained in any fiber (we call such a divisor are horizontal), and
  there exist no common components between $\mcD$ and $\sigma_f^*\mcD$,  and
 
 \item[(ii)] all irreducible components $\Xi$ and $\sigma_f^*\Xi$ are fiber components of $\varphi$ and 
 there exist no common components between $\Xi$ and $\sigma_f^*\Xi$.
 
 \end{enumerate}
 
 \begin{rem}\label{rem:components} { \rm Possible irreducible components of $\Xi$ and $\sigma_f^*\Xi$ can be determined by
  \cite[Remark 3.1]{tokunaga14}. In particular, if singular fibers of $\varphi$ are of types $\I_1, \I_2, \II, \III$ only,
  $\Xi_f = \emptyset$.
 }
 \end{rem}
 
%
%
%
%
%
%
%
%



\section{Proof of Theorem~\ref{thm:comb-1}}
 
We first recall the double cover diagram for $S_{\mcQ, z_o}$. In our case,
$\widehat{\Sigma}_2$ can be successively blown down to $\PP^2$. We denote it by $\overline{q} : \widehat{\Sigma}_2 \to \PP^2$.
We then have the following combined diagram:
\[
\begin{CD}
S' @<{\mu}<< S_{\mcQ, z_o} \\
@V{f'}VV                 @VV{f}V \\
\Sigma_2@<<{q}< \widehat{\Sigma}_2 @>{\overline{q}}>>  \PP^2. 
\end{CD}
\]
Note that $f_{\mcQ, z_o} = \overline{q}\circ f$.

 Theorem~\ref{thm:comb-1} will be proved based on \cite[Theorem 3.2]{tokunaga14} and the following lemma.  


 \begin{lem}\label{lem:key-comb1}{ Let $P_1, P_2, P_3\in E_{\mcQ,z_o}(\CC(t))$ be rational points such that $f_{\mcQ, z_o}(s_{P_i})=\mcL_i$. 

 Then there exists a $D_{2p}$-cover $\pi_p : X_p \to \widehat{\Sigma}_2$ such that
 
 \begin{itemize}
 
   \item $D(X_p/\widehat{\Sigma}_2) = S_{\mcQ, z_o}$ and 
   
   \item $\Delta_{\beta_2(\pi_p)} = \sum_{i=1}^3 s_{P_i} + \sigma_f^*(\sum_{i=1}^3 s_{P_i})$ 
 \end{itemize}
 if and only if there exists a $D_{2p}$-cover $\overline{\pi}_p : \overline {X}_p \to \PP^2$ such that  
 \begin{itemize}
 
  \item $\Delta_{\overline{\pi}_p} = \mcQ + \sum_{i=1}^3 \mcL_i$ 
  and 
  \item the branch locus of $\beta_1(\overline{\pi}_p) = \mcQ$.
  \end{itemize}       
 }
 \end{lem}
 
 \proof Assume that there exists a $D_{2p}$-cover of $\widehat{\Sigma}_2$ described as above. Consider 
 the Stein factorizaiton $\overline{\pi}_p : \overline{X}_p \to \PP^2$ of $\overline{q}\circ \pi_p$. As the branch
 locus of $\overline{\pi}_p$ is $\overline {q}(\Delta_{\pi_p})$, 
 $\Delta_{\overline{\pi}_p} = \mcQ + \sum_{i=1}^3 \mcL_i$ 
  and the covering $\beta_1(\overline{\pi}_p)$ is branched along $\mcQ$. Conversely, suppose that there exists a $D_{2p}$-cover 
 $\overline{\pi}_p : \overline{X}_p \to \PP^2$ satisfying the above condition.
 Take the $\CC(\overline{X}_p)$-normalization $\pi_p : X_p \to \widehat{\Sigma}_2$ of $\widehat{\Sigma}_2$.
 Since $\Delta_{\beta_1(\overline{\pi}_p)} = \mcQ$, $\overline{q}(\Delta_{\beta_1(\pi)_p}) = \mcQ$. Hence we infer 
 that $D(X_p/\widehat{\Sigma}_2)$ is a double cover of $\widehat{\Sigma}_2$ branched along $f(O)$ and the 
 proper transform of $\mcQ$, i.e., $D(X_p/\widehat{\Sigma}_2) = S_{\mcQ, z_o}$ and $\beta_1(\pi_p) = f$. 
 Let $\beta_1(\pi_p) : X_p \to S_{\mcQ, z_o}$ be the $p$-cyclic cover determined by $\pi_p$. By \cite[Remark 3.1, (i)]{tokunaga14},
 any irreducible component of singular fibers can not be contained in the branch locus of $\beta_2(\pi_p)$. Hence
 $\overline{q}\circ f(\Delta_{\beta_2(\pi_p)}) = \sum_{i=1}^3 \mcL_i$ 
 \qed
 
By the above lemma, we will be able to choose sections appropriately in constructing our configurations so that dihedral covers exist or do not exist. The difference in (non-)existence allows us to distinguish our configurations topologically. 
\bigskip

{\bf Proof for Combinatorics 1-(a)} 
$\mcE$: a smooth cubic. Choose $z_o \in \mcE\setminus \{p_1, p_2, p_3\}$ and let $S_{\mcQ, z_o}$ be
the rational elliptic surface as in the Introduction. 
By Section 1,  $E_{\mcQ, z_o}(\CC(t)) \cong
D_4^*\oplus \ZZ/2\ZZ$. We choose generators $P_0, P_1, P_2, P_3$ of the  $D_4^*$ part such that
\[
[\langle P_i, P_j \rangle ] 
= \left [ 
   \begin{array}{cccc}
      2 & 1 & 1& 1 \\
      1 & 1 & \frac 12 & \frac 12 \\
      1 & \frac 12 & 1 & \frac 12 \\
      1 & \frac 12 & \frac 12 & 1 
    \end{array}
   \right ].
\]   
We denote the $2$-torsion point by $P_{\tau}$. Let $s_{P_i}$ be the corresponding section for each $P_i$ $(i = 0, 1, 2, 3, \tau)$.
Since $\langle P_i, P_j\rangle$ is determined by the intersection numbers, $s_{P_i}s_{P_j}$, $s_{P_i}O$, and
$\mathrm{Contr}_v(s_{P_i}),  \mathrm {Contr}_v(s_{P_i}, s_{P_j})$, i.e., 
at which component of each singular fiber $s_{P_i}$ intersects, we have may assume that the following:    

\begin{enumerate}
 
  \item[(i)] By \cite[Theorem 10.8]{shioda90}, we may assume that $s_{P_i}O = 0$ $(i = 0, 1, 2, 3)$.
  \[
  \sum_{v \in  \Red(\varphi)}\mathrm{Contr}_v(s_{P_i}),  
  \sum_{v \in  \Red(\varphi)}\mathrm {Contr}_v(s_{P_i}, s_{P_j}) = 0, 1/2, 1, 3/2, \mbox{\rm or}\,\,  2.
  \]
  \item[(ii)] For $P_i$ $(i = 1, 2, 3)$, $\langle P_i, P_i \rangle = 1$ implies that 
  $\sum_v \mathrm{Contr}_v(s_{P_i}) = 1$. Also for $\{i, j\} \subset \{1, 2, 3\}, i\neq j$, $\langle P_i, \pm P_j \rangle = 1/2$
  implies $s_{P_i}s_{P_j} = 0$, $\sum_v \mathrm{Contr}_v(s_{P_i}, s_{P_j}) = 1/2$.
  
  \item[(iii)] For $P_0$,  $\sum_v \mathrm{Contr}_v(s_{P_0}) = 0$.
  
  \item[(iv)] For  $P_{\tau}$, as $P_{\tau}$ is a torsion, $\langle P_{\tau}, P_{\tau} \rangle = 0$ and we have 
  $\sum_v \mathrm{Contr}_v(s_{P_{\tau}}) = 2$.
  
\end{enumerate}

From the facts as above, we have that 
\[
\begin{array}{cclcc}
\gamma_{\mcQ, z_o} (P_{\tau}) =   (1, 1, 1,1) & &
\gamma_{\mcQ, z_o} (P_0)  =  (0, 0, 0,0) &&
\gamma_{\mcQ, z_o} (P_1)  =  (1, 1, 0,0)  \\
\gamma_{\mcQ, z_o} (P_2)  =  (1, 0, 1,0) & &
\gamma_{\mcQ, z_o} (P_3)  =  (1, 0, 0,1)  \,\, \mbox{\rm or} \,\,(0, 1, 1, 0) & & 
\end{array}
\]
By replacing $P_3$ by $P_3 \dot{+} P_\tau$, if necessary, we may assume that
 $\gamma_{\mcQ, z_o} (P_3)  =  (1, 0, 0,1)$.   Put $P_4:= P_2 \dot {-} P_3\dot{+}P_{\tau}$. Then
$\gamma_{\mcQ, z_o}(P_4) = (1, 1, 0, 0)$. We now label the irreducible components of the singular fibers as in the following figure:
\begin{center}
{\unitlength 0.1in%
\begin{picture}(44.0000,22.0000)(10.0000,-28.0000)%
%
\special{pn 13}%
\special{pa 1084 2550}%
\special{pa 5175 2542}%
\special{fp}%
%
\special{pn 13}%
\special{pa 4390 830}%
\special{pa 4398 861}%
\special{pa 4406 893}%
\special{pa 4413 924}%
\special{pa 4421 956}%
\special{pa 4428 987}%
\special{pa 4436 1018}%
\special{pa 4443 1050}%
\special{pa 4449 1081}%
\special{pa 4456 1112}%
\special{pa 4462 1143}%
\special{pa 4467 1175}%
\special{pa 4477 1237}%
\special{pa 4485 1299}%
\special{pa 4487 1331}%
\special{pa 4489 1362}%
\special{pa 4489 1394}%
\special{pa 4487 1426}%
\special{pa 4484 1459}%
\special{pa 4479 1493}%
\special{pa 4472 1527}%
\special{pa 4462 1560}%
\special{pa 4449 1593}%
\special{pa 4434 1623}%
\special{pa 4415 1650}%
\special{pa 4392 1673}%
\special{pa 4366 1690}%
\special{pa 4335 1700}%
\special{pa 4302 1702}%
\special{pa 4272 1693}%
\special{pa 4247 1673}%
\special{pa 4230 1644}%
\special{pa 4220 1610}%
\special{pa 4219 1573}%
\special{pa 4228 1539}%
\special{pa 4246 1509}%
\special{pa 4270 1488}%
\special{pa 4299 1479}%
\special{pa 4330 1484}%
\special{pa 4361 1499}%
\special{pa 4389 1522}%
\special{pa 4410 1549}%
\special{pa 4424 1579}%
\special{pa 4432 1610}%
\special{pa 4436 1642}%
\special{pa 4438 1674}%
\special{pa 4438 1707}%
\special{pa 4439 1740}%
\special{pa 4439 1805}%
\special{pa 4440 1837}%
\special{pa 4440 2285}%
\special{pa 4439 2316}%
\special{pa 4439 2475}%
\special{pa 4438 2507}%
\special{pa 4438 2692}%
\special{fp}%
%
\special{pn 13}%
\special{pa 4920 830}%
\special{pa 4928 861}%
\special{pa 4936 893}%
\special{pa 4944 924}%
\special{pa 4951 955}%
\special{pa 4959 987}%
\special{pa 4973 1049}%
\special{pa 4980 1081}%
\special{pa 4998 1174}%
\special{pa 5003 1206}%
\special{pa 5008 1237}%
\special{pa 5016 1299}%
\special{pa 5018 1330}%
\special{pa 5020 1362}%
\special{pa 5020 1394}%
\special{pa 5018 1426}%
\special{pa 5015 1459}%
\special{pa 5010 1492}%
\special{pa 5003 1526}%
\special{pa 4993 1560}%
\special{pa 4980 1593}%
\special{pa 4964 1623}%
\special{pa 4946 1650}%
\special{pa 4923 1672}%
\special{pa 4896 1689}%
\special{pa 4866 1700}%
\special{pa 4833 1702}%
\special{pa 4802 1693}%
\special{pa 4778 1673}%
\special{pa 4760 1645}%
\special{pa 4750 1610}%
\special{pa 4749 1574}%
\special{pa 4758 1539}%
\special{pa 4775 1509}%
\special{pa 4800 1488}%
\special{pa 4828 1479}%
\special{pa 4860 1484}%
\special{pa 4890 1499}%
\special{pa 4918 1521}%
\special{pa 4940 1549}%
\special{pa 4954 1578}%
\special{pa 4963 1609}%
\special{pa 4967 1641}%
\special{pa 4969 1673}%
\special{pa 4969 1706}%
\special{pa 4970 1738}%
\special{pa 4971 1771}%
\special{pa 4971 1836}%
\special{pa 4972 1868}%
\special{pa 4972 2220}%
\special{pa 4971 2252}%
\special{pa 4971 2347}%
\special{pa 4970 2379}%
\special{pa 4970 2506}%
\special{pa 4969 2538}%
\special{pa 4969 2692}%
\special{fp}%
\put(51.6000,-11.0000){\makebox(0,0)[lb]{$s_{P_1}$}}%
\put(51.4000,-12.8000){\makebox(0,0)[lb]{$s_{P_2}$}}%
\put(30.7400,-26.5700){\makebox(0,0)[lb]{$O$}}%
\put(11.1000,-24.0000){\makebox(0,0)[lb]{$\Theta_{0,0}$}}%
\put(33.3000,-17.5000){\makebox(0,0)[lb]{$\Theta_{3,1}$}}%
\put(14.4000,-18.6000){\makebox(0,0)[lb]{$\Theta_{0,1}$}}%
\put(32.9000,-23.2000){\makebox(0,0)[lb]{$\Theta_{3,0}$}}%
%
\special{pn 13}%
\special{pa 2010 770}%
\special{pa 1991 799}%
\special{pa 1971 829}%
\special{pa 1933 887}%
\special{pa 1915 916}%
\special{pa 1897 946}%
\special{pa 1879 975}%
\special{pa 1862 1004}%
\special{pa 1846 1034}%
\special{pa 1831 1063}%
\special{pa 1817 1093}%
\special{pa 1804 1122}%
\special{pa 1792 1151}%
\special{pa 1781 1181}%
\special{pa 1772 1210}%
\special{pa 1764 1240}%
\special{pa 1758 1270}%
\special{pa 1753 1299}%
\special{pa 1751 1329}%
\special{pa 1750 1359}%
\special{pa 1750 1388}%
\special{pa 1753 1418}%
\special{pa 1758 1448}%
\special{pa 1764 1478}%
\special{pa 1772 1508}%
\special{pa 1781 1538}%
\special{pa 1791 1568}%
\special{pa 1802 1598}%
\special{pa 1814 1628}%
\special{pa 1826 1659}%
\special{pa 1852 1719}%
\special{pa 1878 1781}%
\special{pa 1890 1811}%
\special{pa 1903 1842}%
\special{pa 1925 1904}%
\special{pa 1935 1935}%
\special{pa 1943 1967}%
\special{pa 1951 1998}%
\special{pa 1958 2029}%
\special{pa 1964 2061}%
\special{pa 1969 2092}%
\special{pa 1977 2156}%
\special{pa 1983 2252}%
\special{pa 1983 2284}%
\special{pa 1984 2316}%
\special{pa 1980 2444}%
\special{pa 1978 2477}%
\special{pa 1976 2509}%
\special{pa 1973 2541}%
\special{pa 1971 2574}%
\special{pa 1965 2638}%
\special{pa 1962 2671}%
\special{pa 1960 2690}%
\special{fp}%
%
\special{pn 13}%
\special{pa 1200 770}%
\special{pa 1218 799}%
\special{pa 1237 829}%
\special{pa 1255 858}%
\special{pa 1272 887}%
\special{pa 1290 917}%
\special{pa 1324 975}%
\special{pa 1340 1005}%
\special{pa 1355 1034}%
\special{pa 1370 1064}%
\special{pa 1384 1093}%
\special{pa 1397 1122}%
\special{pa 1409 1152}%
\special{pa 1420 1181}%
\special{pa 1430 1211}%
\special{pa 1438 1240}%
\special{pa 1445 1270}%
\special{pa 1451 1299}%
\special{pa 1456 1329}%
\special{pa 1459 1359}%
\special{pa 1460 1388}%
\special{pa 1460 1418}%
\special{pa 1458 1448}%
\special{pa 1455 1477}%
\special{pa 1450 1507}%
\special{pa 1444 1537}%
\special{pa 1436 1566}%
\special{pa 1418 1626}%
\special{pa 1407 1656}%
\special{pa 1395 1686}%
\special{pa 1382 1716}%
\special{pa 1369 1745}%
\special{pa 1324 1835}%
\special{pa 1307 1865}%
\special{pa 1291 1895}%
\special{pa 1274 1925}%
\special{pa 1256 1955}%
\special{pa 1239 1985}%
\special{pa 1200 2050}%
\special{fp}%
%
\special{pn 13}%
\special{pa 1380 760}%
\special{pa 1361 789}%
\special{pa 1341 819}%
\special{pa 1303 877}%
\special{pa 1285 906}%
\special{pa 1267 936}%
\special{pa 1249 965}%
\special{pa 1232 994}%
\special{pa 1216 1024}%
\special{pa 1201 1053}%
\special{pa 1187 1083}%
\special{pa 1174 1112}%
\special{pa 1162 1141}%
\special{pa 1151 1171}%
\special{pa 1142 1200}%
\special{pa 1134 1230}%
\special{pa 1128 1260}%
\special{pa 1123 1289}%
\special{pa 1121 1319}%
\special{pa 1120 1349}%
\special{pa 1120 1378}%
\special{pa 1123 1408}%
\special{pa 1128 1438}%
\special{pa 1134 1468}%
\special{pa 1142 1498}%
\special{pa 1151 1528}%
\special{pa 1161 1558}%
\special{pa 1172 1588}%
\special{pa 1184 1618}%
\special{pa 1196 1649}%
\special{pa 1222 1709}%
\special{pa 1248 1771}%
\special{pa 1260 1801}%
\special{pa 1273 1832}%
\special{pa 1295 1894}%
\special{pa 1305 1925}%
\special{pa 1313 1957}%
\special{pa 1321 1988}%
\special{pa 1328 2019}%
\special{pa 1334 2051}%
\special{pa 1339 2082}%
\special{pa 1347 2146}%
\special{pa 1353 2242}%
\special{pa 1353 2274}%
\special{pa 1354 2306}%
\special{pa 1350 2434}%
\special{pa 1348 2467}%
\special{pa 1346 2499}%
\special{pa 1343 2531}%
\special{pa 1341 2564}%
\special{pa 1335 2628}%
\special{pa 1332 2661}%
\special{pa 1330 2680}%
\special{fp}%
%
\special{pn 13}%
\special{pa 1830 740}%
\special{pa 1848 769}%
\special{pa 1867 799}%
\special{pa 1885 828}%
\special{pa 1902 857}%
\special{pa 1920 887}%
\special{pa 1954 945}%
\special{pa 1970 975}%
\special{pa 1985 1004}%
\special{pa 2000 1034}%
\special{pa 2014 1063}%
\special{pa 2027 1092}%
\special{pa 2039 1122}%
\special{pa 2050 1151}%
\special{pa 2060 1181}%
\special{pa 2068 1210}%
\special{pa 2075 1240}%
\special{pa 2081 1269}%
\special{pa 2086 1299}%
\special{pa 2089 1329}%
\special{pa 2090 1358}%
\special{pa 2090 1388}%
\special{pa 2088 1418}%
\special{pa 2085 1447}%
\special{pa 2080 1477}%
\special{pa 2074 1507}%
\special{pa 2066 1536}%
\special{pa 2048 1596}%
\special{pa 2037 1626}%
\special{pa 2025 1656}%
\special{pa 2012 1686}%
\special{pa 1999 1715}%
\special{pa 1954 1805}%
\special{pa 1937 1835}%
\special{pa 1921 1865}%
\special{pa 1904 1895}%
\special{pa 1886 1925}%
\special{pa 1869 1955}%
\special{pa 1830 2020}%
\special{fp}%
%
\special{pn 13}%
\special{pa 2670 790}%
\special{pa 2651 819}%
\special{pa 2631 849}%
\special{pa 2593 907}%
\special{pa 2575 936}%
\special{pa 2557 966}%
\special{pa 2539 995}%
\special{pa 2522 1024}%
\special{pa 2506 1054}%
\special{pa 2491 1083}%
\special{pa 2477 1113}%
\special{pa 2464 1142}%
\special{pa 2452 1171}%
\special{pa 2441 1201}%
\special{pa 2432 1230}%
\special{pa 2424 1260}%
\special{pa 2418 1290}%
\special{pa 2413 1319}%
\special{pa 2411 1349}%
\special{pa 2410 1379}%
\special{pa 2410 1408}%
\special{pa 2413 1438}%
\special{pa 2418 1468}%
\special{pa 2424 1498}%
\special{pa 2432 1528}%
\special{pa 2441 1558}%
\special{pa 2451 1588}%
\special{pa 2462 1618}%
\special{pa 2474 1648}%
\special{pa 2486 1679}%
\special{pa 2512 1739}%
\special{pa 2538 1801}%
\special{pa 2550 1831}%
\special{pa 2563 1862}%
\special{pa 2585 1924}%
\special{pa 2595 1955}%
\special{pa 2603 1987}%
\special{pa 2611 2018}%
\special{pa 2618 2049}%
\special{pa 2624 2081}%
\special{pa 2629 2112}%
\special{pa 2637 2176}%
\special{pa 2643 2272}%
\special{pa 2643 2304}%
\special{pa 2644 2336}%
\special{pa 2640 2464}%
\special{pa 2638 2497}%
\special{pa 2636 2529}%
\special{pa 2633 2561}%
\special{pa 2631 2594}%
\special{pa 2625 2658}%
\special{pa 2622 2691}%
\special{pa 2620 2710}%
\special{fp}%
%
\special{pn 13}%
\special{pa 2500 750}%
\special{pa 2518 779}%
\special{pa 2537 809}%
\special{pa 2555 838}%
\special{pa 2572 867}%
\special{pa 2590 897}%
\special{pa 2624 955}%
\special{pa 2640 985}%
\special{pa 2655 1014}%
\special{pa 2670 1044}%
\special{pa 2684 1073}%
\special{pa 2697 1102}%
\special{pa 2709 1132}%
\special{pa 2720 1161}%
\special{pa 2730 1191}%
\special{pa 2738 1220}%
\special{pa 2745 1250}%
\special{pa 2751 1279}%
\special{pa 2756 1309}%
\special{pa 2759 1339}%
\special{pa 2760 1368}%
\special{pa 2760 1398}%
\special{pa 2758 1428}%
\special{pa 2755 1457}%
\special{pa 2750 1487}%
\special{pa 2744 1517}%
\special{pa 2736 1546}%
\special{pa 2718 1606}%
\special{pa 2707 1636}%
\special{pa 2695 1666}%
\special{pa 2682 1696}%
\special{pa 2669 1725}%
\special{pa 2624 1815}%
\special{pa 2607 1845}%
\special{pa 2591 1875}%
\special{pa 2574 1905}%
\special{pa 2556 1935}%
\special{pa 2539 1965}%
\special{pa 2500 2030}%
\special{fp}%
%
\special{pn 13}%
\special{pa 3290 790}%
\special{pa 3271 819}%
\special{pa 3251 849}%
\special{pa 3213 907}%
\special{pa 3195 936}%
\special{pa 3177 966}%
\special{pa 3159 995}%
\special{pa 3142 1024}%
\special{pa 3126 1054}%
\special{pa 3111 1083}%
\special{pa 3097 1113}%
\special{pa 3084 1142}%
\special{pa 3072 1171}%
\special{pa 3061 1201}%
\special{pa 3052 1230}%
\special{pa 3044 1260}%
\special{pa 3038 1290}%
\special{pa 3033 1319}%
\special{pa 3031 1349}%
\special{pa 3030 1379}%
\special{pa 3030 1408}%
\special{pa 3033 1438}%
\special{pa 3038 1468}%
\special{pa 3044 1498}%
\special{pa 3052 1528}%
\special{pa 3061 1558}%
\special{pa 3071 1588}%
\special{pa 3082 1618}%
\special{pa 3094 1648}%
\special{pa 3106 1679}%
\special{pa 3132 1739}%
\special{pa 3158 1801}%
\special{pa 3170 1831}%
\special{pa 3183 1862}%
\special{pa 3205 1924}%
\special{pa 3215 1955}%
\special{pa 3223 1987}%
\special{pa 3231 2018}%
\special{pa 3238 2049}%
\special{pa 3244 2081}%
\special{pa 3249 2112}%
\special{pa 3257 2176}%
\special{pa 3263 2272}%
\special{pa 3263 2304}%
\special{pa 3264 2336}%
\special{pa 3260 2464}%
\special{pa 3258 2497}%
\special{pa 3256 2529}%
\special{pa 3253 2561}%
\special{pa 3251 2594}%
\special{pa 3245 2658}%
\special{pa 3242 2691}%
\special{pa 3240 2710}%
\special{fp}%
%
\special{pn 13}%
\special{pa 3120 760}%
\special{pa 3138 789}%
\special{pa 3157 819}%
\special{pa 3175 848}%
\special{pa 3192 877}%
\special{pa 3210 907}%
\special{pa 3244 965}%
\special{pa 3260 995}%
\special{pa 3275 1024}%
\special{pa 3290 1054}%
\special{pa 3304 1083}%
\special{pa 3317 1112}%
\special{pa 3329 1142}%
\special{pa 3340 1171}%
\special{pa 3350 1201}%
\special{pa 3358 1230}%
\special{pa 3365 1260}%
\special{pa 3371 1289}%
\special{pa 3376 1319}%
\special{pa 3379 1349}%
\special{pa 3380 1378}%
\special{pa 3380 1408}%
\special{pa 3378 1438}%
\special{pa 3375 1467}%
\special{pa 3370 1497}%
\special{pa 3364 1527}%
\special{pa 3356 1556}%
\special{pa 3338 1616}%
\special{pa 3327 1646}%
\special{pa 3315 1676}%
\special{pa 3302 1706}%
\special{pa 3289 1735}%
\special{pa 3244 1825}%
\special{pa 3227 1855}%
\special{pa 3211 1885}%
\special{pa 3194 1915}%
\special{pa 3176 1945}%
\special{pa 3159 1975}%
\special{pa 3120 2040}%
\special{fp}%
\put(20.1000,-19.3000){\makebox(0,0)[lb]{$\Theta_{1,1}$}}%
\put(19.9000,-22.5000){\makebox(0,0)[lb]{$\Theta_{1,0}$}}%
\put(26.6000,-24.9000){\makebox(0,0)[lb]{$\Theta_{2,0}$}}%
\put(26.4000,-19.9000){\makebox(0,0)[lb]{$\Theta_{2,1}$}}%
\put(50.9000,-9.2000){\makebox(0,0)[lb]{$s_{P_0}$}}%
\put(51.2000,-17.3000){\makebox(0,0)[lb]{$s_{P_{\tau}}$}}%
%
\special{pn 13}%
\special{pa 1000 600}%
\special{pa 1000 2800}%
\special{fp}%
\special{pa 5400 2800}%
\special{pa 5400 2800}%
\special{fp}%
\special{pa 1000 2800}%
\special{pa 5400 2800}%
\special{fp}%
\special{pa 5400 2800}%
\special{pa 5400 600}%
\special{fp}%
\special{pa 5400 600}%
\special{pa 1000 600}%
\special{fp}%
%
\special{pn 13}%
\special{pa 1100 1010}%
\special{pa 1300 1010}%
\special{fp}%
\special{pa 1430 1010}%
\special{pa 1920 1010}%
\special{fp}%
\special{pa 2080 1010}%
\special{pa 2590 1010}%
\special{fp}%
\special{pa 2750 1010}%
\special{pa 3220 1010}%
\special{fp}%
\special{pa 3370 1010}%
\special{pa 5170 1010}%
\special{fp}%
%
\special{pn 13}%
\special{pa 1255 1130}%
\special{pa 1715 1130}%
\special{fp}%
\special{pa 1885 1130}%
\special{pa 2635 1130}%
\special{fp}%
\special{pa 2805 1130}%
\special{pa 3245 1130}%
\special{fp}%
\special{pa 3435 1130}%
\special{pa 5205 1130}%
\special{fp}%
%
\special{pn 13}%
\special{pa 1220 1280}%
\special{pa 2020 1280}%
\special{fp}%
\special{pa 2200 1280}%
\special{pa 2330 1280}%
\special{fp}%
\special{pa 2510 1280}%
\special{pa 3310 1280}%
\special{fp}%
\special{pa 3460 1280}%
\special{pa 5220 1280}%
\special{fp}%
%
\special{pn 13}%
\special{pa 1250 1430}%
\special{pa 2010 1430}%
\special{fp}%
\special{pa 2180 1430}%
\special{pa 2660 1430}%
\special{fp}%
\special{pa 2820 1430}%
\special{pa 2960 1430}%
\special{fp}%
\special{pa 3120 1430}%
\special{pa 5250 1430}%
\special{fp}%
%
\special{pn 13}%
\special{pa 1285 1590}%
\special{pa 1705 1590}%
\special{fp}%
\special{pa 1875 1590}%
\special{pa 2365 1590}%
\special{fp}%
\special{pa 2545 1590}%
\special{pa 2985 1590}%
\special{fp}%
\special{pa 3155 1590}%
\special{pa 4345 1590}%
\special{fp}%
\special{pa 4535 1590}%
\special{pa 4865 1590}%
\special{fp}%
\special{pa 5065 1590}%
\special{pa 5205 1590}%
\special{fp}%
\put(51.7000,-15.3000){\makebox(0,0)[lb]{$s_{P_3}$}}%
%
\special{pn 4}%
\special{sh 1}%
\special{ar 4600 2100 10 10 0 6.2831853}%
\special{sh 1}%
\special{ar 4800 2100 10 10 0 6.2831853}%
\special{sh 1}%
\special{ar 4700 2100 10 10 0 6.2831853}%
\special{sh 1}%
\special{ar 4700 2100 10 10 0 6.2831853}%
\end{picture}}%

Figure 1
\end{center}
%
(Note that we label singular fibers of type $\III$ similarly, if they exist.) 
We now blow down  smooth rational curves $f(\Theta_{0, 0}), f(O), f(\Theta_{1, 1}),
f(\Theta_{2, 1}),  f(\Theta_{3, 1})$  in this order. The resulting surface is $\PP^2$ and  this is the morphism
$\overline q$ in this case.
Note that $f_{\mcQ, z_o} = \overline{q}\circ f$ and $\overline {q}\circ f(O\cup \Theta_{0,0}) = z_o$. 


\begin{lem}\label{lem:image-1}{
\begin{enumerate}
 \item[(i)]  The image of the fixed locus of $[-1]_{\varphi}$ is a smooth cubic $\mcE$ and a transversal line $\mcL_o$ to $\mcE$. Moreover, $\overline {q}\circ f(s_{P_{\tau}}) = \mcL_o$.
 
 \item[(ii)]  $\{f_{\mcQ, z_o}(\Theta_{1,1}), f_{\mcQ, z_o}(\Theta_{2,1}), f_{\mcQ, z_o}(\Theta_{3,1})\} =
 \mcE\cap \mcL_o$. We denote $p_i = f_{\mcQ, z_o}(\Theta_{i,1})$.
 
\item[(iii)]  Put $\mcL_i := f_{\mcQ, z_o}(s_{P_i})$. Then $f_{\mcQ, z_o}(s_{P_i})$ $(i =  1, 2, 3, 4)$ passes through $p_i$ ($i = 1, 2, 3$), respectively.
 Also  $\mcL_4$ passes through $p_1$.

 \item[(iv)] Furthermore, 
 (a) $\mcL_i$ is tangent to $\mcE$ at a point other than $p_i$ ($p_1$ for $\mcL_4$)   or 
 $p_i$ (resp. $p_1$) is an inflection point of $\mcE$ and $\mcL_i$ (resp. $\mcL_4$) is
 an inflectional tangent line.
 
 \end{enumerate}
 }
 \end{lem}
 
 \proof 
 The statements (i) and (ii) are immediate from our construction of $S_{\mcQ, z_o}$.
%
%
 We show that the statements (iii) and (iv) hold for $\mcL_1= f_{\mcQ, z_o}(s_{P_1})$ only, as our proof for the remaining sections are the same.
 Since $\gamma_{\mcQ, z_o}(P_1) = (1, 1, 0, 0)$, $s_{P_1}\Theta_{2,1} = 1$.
 This shows that  $f_{\mcQ, z_o}(\Theta_{2,1}) \in \mcL_1$. 
 We now go on to (iv). 
 Since $\gamma_{\mcQ, z_o}(P_1) = (1, 1, 0, 0)$, $s_{P_1}\Theta_{1,0} = 0$ and $z_o \not\in \mcL_1$.
 As general fiber
 $F$ fo $\varphi_{\mcQ, z_o}$ is mapped to a line $l$ through $z_o$ and $s_{P_1}F = 1$,
 $z_o \not\in f_{\mcQ, z_o}(s_{P_1})$ implies that $\mcL_1\cap l$ consists of only one point.  This shows that
 $\mcL_1$ is a line. If $\mcL_1$ is not the line described in (iv) , the closure
 of $f_{\mcQ, z_o}^{-1}(\mcL_1 \setminus ( \mcL_1\cap \mcE))$ is irreducible. on the other hand, it must contains
 both $s_{P_1}$ and $[-1]_{\varphi_{\mcQ, z_o}}(s_{P_1})$, which is a contradiction. \qed
  
\begin{remark}\label{rem:obs-1} Conversely, from the proof of Lemma \ref{lem:image-1} we observe that:

 \begin{enumerate}
 \item[(i)] Any tangent line to $\mcE$ through $p_i$ gives rise a point $P \in E_{\mcQ, z_o}(\CC(t))$ with $\langle P, P \rangle = 1$.
 
 \item[(ii)] Any contact conic to $\mcQ$ through $z_o$ gives rise a point $P \in E_{\mcQ, z_o}(\CC(t))$ with $\langle P, P \rangle = 2$.
 
 \end{enumerate}
\end{remark}
  
Now let
\[
\mcB^1=\mcE +\mcL_0+\sum_{i=1}^3 \mcL_i, 
\mcB^2=\mcE +\mcL_0+\sum_{i=2}^4 \mcL_i.
\]
Then by  \cite[Theorem 3.2]{tokunaga14} and  Lemma~\ref{lem:key-comb1}, a $D_{2p}$-cover branched at $2(\mcE+\mcL_0)+p(\sum_{i=2}^4 \mcL_i)$ exists, but does not exist for $2(\mcE +\mcL_0)+p(\sum_{i=1}^3 \mcL_i)$. Hence, $(\mcB^1, \mcB^2)$  is a Zariski pair
if their combinatorics are the same.
 \bigskip

 \bigskip

{\bf Proof for Combinatorics 1-(b)}  $\mcE$: a nodal cubic. Choose $z_o \in \mcE\setminus \{p_1, p_2, p_3\}$ and let $S_{\mcQ, z_o}$ be
the rational elliptic surface as in the Introduction. 
By Section 1, $E_{\mcQ, z_o}(\CC(t)) \cong
(A_1^*)^{\oplus 3}\oplus \ZZ/2\ZZ$. We choose generators $P_1, P_2, P_3$ of the  $(A_1^*)^{\oplus 3}$-part
such that
\[
[\langle P_i, P_j\rangle]=\left[ \begin{array}{ccc} \frac{1}{2}  & 0 & 0 \\ 0 & \frac{1}{2}  & 0 \\ 0 & 0 &\frac{1}{2} \end{array}\right]
\]
We denote the 2-torsion point by $P_{\tau}$. Let $s_{P_i}$ be the corresponding section for each $P_i$ ($i=1,2,3,\tau$)

\begin{enumerate}
\item[(i)] By \cite[Theorem 10.8]{shioda90}, we may assume that $s_{P_i}O=0$, $(i=0,1, 2,\tau)$ and
\[
\sum_{v \in  \Red(\varphi)}\mathrm{Contr}_v(s_{P_i}),  
  \sum_{v \in  \Red(\varphi)}\mathrm {Contr}_v(s_{P_i}, s_{P_j}) = 0, 1/2, 1, 3/2, \mbox{\rm or}\,\,  2, 5/2.
\]
\item[(ii)] For $P_i$ ($i=1,2,3$), $\langle P_i, P_i\rangle= 1/2$ implies that $\sum_v\mathrm{Contr}_v(s_{P_i})=3/2$. Also for $\left\{i,j\right\}\subset\{1,2,3\}$, $i\not=j$, $\langle P_i, P_j\rangle=0$ implies $s_{P_i}s_{P_j}=0$, $\sum_v\mathrm{Contr}_v(s_{P_i}, s_{P_j})=1$
\item[(iii)] For $P_\tau$, as $P_\tau$ is a torsion section, $ \langle P_\tau, P_\tau\rangle=0$ and we have $\sum_v\mathrm{Contr}_v(s_{P_i})=2$.
\item[(iv)] Furthermore, $\langle P_i, P_\tau\rangle=0$ implies  $\sum_v\mathrm{Contr}_v(s_{P_i}, s_{P_\tau})=1$.
\end{enumerate}

From the above facts, we can assume that
\[
\begin{array}{ccl}
\gamma_{\mcQ, z_o} (P_{\tau})  =  (1, 1, 1, 1 , 0) & &
\gamma_{\mcQ, z_o} (P_1)  =  (1, 1, 0, 0, 1) \\
\gamma_{\mcQ, z_o} (P_2)  =  (1, 0, 1, 0, 1) & &
\gamma_{\mcQ, z_o} (P_3)  =  (1, 0, 0,1, 1)  \,\, \mbox{\rm or} \,\,(0, 1, 1, 0, 1)
\end{array}
\]
By replacing $P_3$  by $P_3\dot+P_\tau$  if necessary, we may assume that $\gamma_{\mcQ, z_o} (P_3)  =  (1, 0, 0,1, 1).$ We now label the irreducible components of the singular fibers as in the following figure. (cf.
\cite[No. 24, p. 90]{tokunaga12}. Note that we use different labelings.). 

%


Now consider $P_{i,j}^\pm=P_i\dot\pm P_j\dot+P_\tau$ ($\{i,j\}\subset\{1,2,3\}$, $i\not=j$). Then, since $\gamma_{\mcQ, z_o}$ is a group homomorphism, $s_{P_{i,j}^\pm}\Theta_{1,1}=s_{P_{i,j}^\pm}\Theta_{k,1}=1$, $k\not=i,j$. Also, since $\sum_v{\rm Contr}(s_{P_{i,j}^\pm})=1$, and $\langle P_{i,j}^+, P_{i,j}^+\rangle=\langle P_{i,j}^-, P_{i,j}^-\rangle=1$, we have $s_{P_{i,j}^\pm}O=0$.

Let $f: S_{\mcQ, z_o} \to \widehat{\Sigma}_2$ be the one in the
double cover diagram.
We now blow down  smooth rational curves $f(\Theta_{0, 0})$, $f(O)$, $f(\Theta_{1, 1})$,
$f(\Theta_{2, 1})$, $f(\Theta_{3, 1})$, $f(\Theta_{4, 1})$  in this order. The resulting surface is $\PP^2$ and this morphism coincides with $\overline q$ as in the previous case.

\begin{center}
{\unitlength 0.1in%
\begin{picture}(43.8000,23.9000)(9.8000,-28.1000)%
%
\special{pn 13}%
\special{pa 980 420}%
\special{pa 5360 420}%
\special{pa 5360 2810}%
\special{pa 980 2810}%
\special{pa 980 420}%
\special{pa 5360 420}%
\special{fp}%
%
\special{pn 13}%
\special{pa 1084 2550}%
\special{pa 5175 2542}%
\special{fp}%
%
\special{pn 13}%
\special{pa 4390 830}%
\special{pa 4398 861}%
\special{pa 4406 893}%
\special{pa 4413 924}%
\special{pa 4421 956}%
\special{pa 4428 987}%
\special{pa 4436 1018}%
\special{pa 4443 1050}%
\special{pa 4449 1081}%
\special{pa 4456 1112}%
\special{pa 4462 1143}%
\special{pa 4467 1175}%
\special{pa 4477 1237}%
\special{pa 4485 1299}%
\special{pa 4487 1331}%
\special{pa 4489 1362}%
\special{pa 4489 1394}%
\special{pa 4487 1426}%
\special{pa 4484 1459}%
\special{pa 4479 1493}%
\special{pa 4472 1527}%
\special{pa 4462 1560}%
\special{pa 4449 1593}%
\special{pa 4434 1623}%
\special{pa 4415 1650}%
\special{pa 4392 1673}%
\special{pa 4366 1690}%
\special{pa 4335 1700}%
\special{pa 4302 1702}%
\special{pa 4272 1693}%
\special{pa 4247 1673}%
\special{pa 4230 1644}%
\special{pa 4220 1610}%
\special{pa 4219 1573}%
\special{pa 4228 1539}%
\special{pa 4246 1509}%
\special{pa 4270 1488}%
\special{pa 4299 1479}%
\special{pa 4330 1484}%
\special{pa 4361 1499}%
\special{pa 4389 1522}%
\special{pa 4410 1549}%
\special{pa 4424 1579}%
\special{pa 4432 1610}%
\special{pa 4436 1642}%
\special{pa 4438 1674}%
\special{pa 4438 1707}%
\special{pa 4439 1740}%
\special{pa 4439 1805}%
\special{pa 4440 1837}%
\special{pa 4440 2285}%
\special{pa 4439 2316}%
\special{pa 4439 2475}%
\special{pa 4438 2507}%
\special{pa 4438 2692}%
\special{fp}%
%
\special{pn 13}%
\special{pa 4920 830}%
\special{pa 4928 861}%
\special{pa 4936 893}%
\special{pa 4944 924}%
\special{pa 4951 955}%
\special{pa 4959 987}%
\special{pa 4973 1049}%
\special{pa 4980 1081}%
\special{pa 4998 1174}%
\special{pa 5003 1206}%
\special{pa 5008 1237}%
\special{pa 5016 1299}%
\special{pa 5018 1330}%
\special{pa 5020 1362}%
\special{pa 5020 1394}%
\special{pa 5018 1426}%
\special{pa 5015 1459}%
\special{pa 5010 1492}%
\special{pa 5003 1526}%
\special{pa 4993 1560}%
\special{pa 4980 1593}%
\special{pa 4964 1623}%
\special{pa 4946 1650}%
\special{pa 4923 1672}%
\special{pa 4896 1689}%
\special{pa 4866 1700}%
\special{pa 4833 1702}%
\special{pa 4802 1693}%
\special{pa 4778 1673}%
\special{pa 4760 1645}%
\special{pa 4750 1610}%
\special{pa 4749 1574}%
\special{pa 4758 1539}%
\special{pa 4775 1509}%
\special{pa 4800 1488}%
\special{pa 4828 1479}%
\special{pa 4860 1484}%
\special{pa 4890 1499}%
\special{pa 4918 1521}%
\special{pa 4940 1549}%
\special{pa 4954 1578}%
\special{pa 4963 1609}%
\special{pa 4967 1641}%
\special{pa 4969 1673}%
\special{pa 4969 1706}%
\special{pa 4970 1738}%
\special{pa 4971 1771}%
\special{pa 4971 1836}%
\special{pa 4972 1868}%
\special{pa 4972 2220}%
\special{pa 4971 2252}%
\special{pa 4971 2347}%
\special{pa 4970 2379}%
\special{pa 4970 2506}%
\special{pa 4969 2538}%
\special{pa 4969 2692}%
\special{fp}%
\put(51.5500,-12.2500){\makebox(0,0)[lb]{$s_{P_2}$}}%
\put(50.4500,-14.8000){\makebox(0,0)[lb]{$s_{P_3}$}}%
\put(30.7400,-26.5700){\makebox(0,0)[lb]{$O$}}%
\put(11.1000,-24.0000){\makebox(0,0)[lb]{$\Theta_{0,0}$}}%
\put(33.3000,-17.5000){\makebox(0,0)[lb]{$\Theta_{3,1}$}}%
\put(14.4000,-18.6000){\makebox(0,0)[lb]{$\Theta_{0,1}$}}%
\put(32.9000,-23.2000){\makebox(0,0)[lb]{$\Theta_{3,0}$}}%
%
\special{pn 13}%
\special{pa 2010 770}%
\special{pa 1991 799}%
\special{pa 1971 829}%
\special{pa 1933 887}%
\special{pa 1915 916}%
\special{pa 1897 946}%
\special{pa 1879 975}%
\special{pa 1862 1004}%
\special{pa 1846 1034}%
\special{pa 1831 1063}%
\special{pa 1817 1093}%
\special{pa 1804 1122}%
\special{pa 1792 1151}%
\special{pa 1781 1181}%
\special{pa 1772 1210}%
\special{pa 1764 1240}%
\special{pa 1758 1270}%
\special{pa 1753 1299}%
\special{pa 1751 1329}%
\special{pa 1750 1359}%
\special{pa 1750 1388}%
\special{pa 1753 1418}%
\special{pa 1758 1448}%
\special{pa 1764 1478}%
\special{pa 1772 1508}%
\special{pa 1781 1538}%
\special{pa 1791 1568}%
\special{pa 1802 1598}%
\special{pa 1814 1628}%
\special{pa 1826 1659}%
\special{pa 1852 1719}%
\special{pa 1878 1781}%
\special{pa 1890 1811}%
\special{pa 1903 1842}%
\special{pa 1925 1904}%
\special{pa 1935 1935}%
\special{pa 1943 1967}%
\special{pa 1951 1998}%
\special{pa 1958 2029}%
\special{pa 1964 2061}%
\special{pa 1969 2092}%
\special{pa 1977 2156}%
\special{pa 1983 2252}%
\special{pa 1983 2284}%
\special{pa 1984 2316}%
\special{pa 1980 2444}%
\special{pa 1978 2477}%
\special{pa 1976 2509}%
\special{pa 1973 2541}%
\special{pa 1971 2574}%
\special{pa 1965 2638}%
\special{pa 1962 2671}%
\special{pa 1960 2690}%
\special{fp}%
%
\special{pn 13}%
\special{pa 1200 770}%
\special{pa 1218 799}%
\special{pa 1237 829}%
\special{pa 1255 858}%
\special{pa 1272 887}%
\special{pa 1290 917}%
\special{pa 1324 975}%
\special{pa 1340 1005}%
\special{pa 1355 1034}%
\special{pa 1370 1064}%
\special{pa 1384 1093}%
\special{pa 1397 1122}%
\special{pa 1409 1152}%
\special{pa 1420 1181}%
\special{pa 1430 1211}%
\special{pa 1438 1240}%
\special{pa 1445 1270}%
\special{pa 1451 1299}%
\special{pa 1456 1329}%
\special{pa 1459 1359}%
\special{pa 1460 1388}%
\special{pa 1460 1418}%
\special{pa 1458 1448}%
\special{pa 1455 1477}%
\special{pa 1450 1507}%
\special{pa 1444 1537}%
\special{pa 1436 1566}%
\special{pa 1418 1626}%
\special{pa 1407 1656}%
\special{pa 1395 1686}%
\special{pa 1382 1716}%
\special{pa 1369 1745}%
\special{pa 1324 1835}%
\special{pa 1307 1865}%
\special{pa 1291 1895}%
\special{pa 1274 1925}%
\special{pa 1256 1955}%
\special{pa 1239 1985}%
\special{pa 1200 2050}%
\special{fp}%
%
\special{pn 13}%
\special{pa 1380 760}%
\special{pa 1361 789}%
\special{pa 1341 819}%
\special{pa 1303 877}%
\special{pa 1285 906}%
\special{pa 1267 936}%
\special{pa 1249 965}%
\special{pa 1232 994}%
\special{pa 1216 1024}%
\special{pa 1201 1053}%
\special{pa 1187 1083}%
\special{pa 1174 1112}%
\special{pa 1162 1141}%
\special{pa 1151 1171}%
\special{pa 1142 1200}%
\special{pa 1134 1230}%
\special{pa 1128 1260}%
\special{pa 1123 1289}%
\special{pa 1121 1319}%
\special{pa 1120 1349}%
\special{pa 1120 1378}%
\special{pa 1123 1408}%
\special{pa 1128 1438}%
\special{pa 1134 1468}%
\special{pa 1142 1498}%
\special{pa 1151 1528}%
\special{pa 1161 1558}%
\special{pa 1172 1588}%
\special{pa 1184 1618}%
\special{pa 1196 1649}%
\special{pa 1222 1709}%
\special{pa 1248 1771}%
\special{pa 1260 1801}%
\special{pa 1273 1832}%
\special{pa 1295 1894}%
\special{pa 1305 1925}%
\special{pa 1313 1957}%
\special{pa 1321 1988}%
\special{pa 1328 2019}%
\special{pa 1334 2051}%
\special{pa 1339 2082}%
\special{pa 1347 2146}%
\special{pa 1353 2242}%
\special{pa 1353 2274}%
\special{pa 1354 2306}%
\special{pa 1350 2434}%
\special{pa 1348 2467}%
\special{pa 1346 2499}%
\special{pa 1343 2531}%
\special{pa 1341 2564}%
\special{pa 1335 2628}%
\special{pa 1332 2661}%
\special{pa 1330 2680}%
\special{fp}%
%
\special{pn 13}%
\special{pa 1830 740}%
\special{pa 1848 769}%
\special{pa 1867 799}%
\special{pa 1885 828}%
\special{pa 1902 857}%
\special{pa 1920 887}%
\special{pa 1954 945}%
\special{pa 1970 975}%
\special{pa 1985 1004}%
\special{pa 2000 1034}%
\special{pa 2014 1063}%
\special{pa 2027 1092}%
\special{pa 2039 1122}%
\special{pa 2050 1151}%
\special{pa 2060 1181}%
\special{pa 2068 1210}%
\special{pa 2075 1240}%
\special{pa 2081 1269}%
\special{pa 2086 1299}%
\special{pa 2089 1329}%
\special{pa 2090 1358}%
\special{pa 2090 1388}%
\special{pa 2088 1418}%
\special{pa 2085 1447}%
\special{pa 2080 1477}%
\special{pa 2074 1507}%
\special{pa 2066 1536}%
\special{pa 2048 1596}%
\special{pa 2037 1626}%
\special{pa 2025 1656}%
\special{pa 2012 1686}%
\special{pa 1999 1715}%
\special{pa 1954 1805}%
\special{pa 1937 1835}%
\special{pa 1921 1865}%
\special{pa 1904 1895}%
\special{pa 1886 1925}%
\special{pa 1869 1955}%
\special{pa 1830 2020}%
\special{fp}%
%
\special{pn 13}%
\special{pa 2670 790}%
\special{pa 2651 819}%
\special{pa 2631 849}%
\special{pa 2593 907}%
\special{pa 2575 936}%
\special{pa 2557 966}%
\special{pa 2539 995}%
\special{pa 2522 1024}%
\special{pa 2506 1054}%
\special{pa 2491 1083}%
\special{pa 2477 1113}%
\special{pa 2464 1142}%
\special{pa 2452 1171}%
\special{pa 2441 1201}%
\special{pa 2432 1230}%
\special{pa 2424 1260}%
\special{pa 2418 1290}%
\special{pa 2413 1319}%
\special{pa 2411 1349}%
\special{pa 2410 1379}%
\special{pa 2410 1408}%
\special{pa 2413 1438}%
\special{pa 2418 1468}%
\special{pa 2424 1498}%
\special{pa 2432 1528}%
\special{pa 2441 1558}%
\special{pa 2451 1588}%
\special{pa 2462 1618}%
\special{pa 2474 1648}%
\special{pa 2486 1679}%
\special{pa 2512 1739}%
\special{pa 2538 1801}%
\special{pa 2550 1831}%
\special{pa 2563 1862}%
\special{pa 2585 1924}%
\special{pa 2595 1955}%
\special{pa 2603 1987}%
\special{pa 2611 2018}%
\special{pa 2618 2049}%
\special{pa 2624 2081}%
\special{pa 2629 2112}%
\special{pa 2637 2176}%
\special{pa 2643 2272}%
\special{pa 2643 2304}%
\special{pa 2644 2336}%
\special{pa 2640 2464}%
\special{pa 2638 2497}%
\special{pa 2636 2529}%
\special{pa 2633 2561}%
\special{pa 2631 2594}%
\special{pa 2625 2658}%
\special{pa 2622 2691}%
\special{pa 2620 2710}%
\special{fp}%
%
\special{pn 13}%
\special{pa 2500 750}%
\special{pa 2518 779}%
\special{pa 2537 809}%
\special{pa 2555 838}%
\special{pa 2572 867}%
\special{pa 2590 897}%
\special{pa 2624 955}%
\special{pa 2640 985}%
\special{pa 2655 1014}%
\special{pa 2670 1044}%
\special{pa 2684 1073}%
\special{pa 2697 1102}%
\special{pa 2709 1132}%
\special{pa 2720 1161}%
\special{pa 2730 1191}%
\special{pa 2738 1220}%
\special{pa 2745 1250}%
\special{pa 2751 1279}%
\special{pa 2756 1309}%
\special{pa 2759 1339}%
\special{pa 2760 1368}%
\special{pa 2760 1398}%
\special{pa 2758 1428}%
\special{pa 2755 1457}%
\special{pa 2750 1487}%
\special{pa 2744 1517}%
\special{pa 2736 1546}%
\special{pa 2718 1606}%
\special{pa 2707 1636}%
\special{pa 2695 1666}%
\special{pa 2682 1696}%
\special{pa 2669 1725}%
\special{pa 2624 1815}%
\special{pa 2607 1845}%
\special{pa 2591 1875}%
\special{pa 2574 1905}%
\special{pa 2556 1935}%
\special{pa 2539 1965}%
\special{pa 2500 2030}%
\special{fp}%
%
\special{pn 13}%
\special{pa 3290 790}%
\special{pa 3271 819}%
\special{pa 3251 849}%
\special{pa 3213 907}%
\special{pa 3195 936}%
\special{pa 3177 966}%
\special{pa 3159 995}%
\special{pa 3142 1024}%
\special{pa 3126 1054}%
\special{pa 3111 1083}%
\special{pa 3097 1113}%
\special{pa 3084 1142}%
\special{pa 3072 1171}%
\special{pa 3061 1201}%
\special{pa 3052 1230}%
\special{pa 3044 1260}%
\special{pa 3038 1290}%
\special{pa 3033 1319}%
\special{pa 3031 1349}%
\special{pa 3030 1379}%
\special{pa 3030 1408}%
\special{pa 3033 1438}%
\special{pa 3038 1468}%
\special{pa 3044 1498}%
\special{pa 3052 1528}%
\special{pa 3061 1558}%
\special{pa 3071 1588}%
\special{pa 3082 1618}%
\special{pa 3094 1648}%
\special{pa 3106 1679}%
\special{pa 3132 1739}%
\special{pa 3158 1801}%
\special{pa 3170 1831}%
\special{pa 3183 1862}%
\special{pa 3205 1924}%
\special{pa 3215 1955}%
\special{pa 3223 1987}%
\special{pa 3231 2018}%
\special{pa 3238 2049}%
\special{pa 3244 2081}%
\special{pa 3249 2112}%
\special{pa 3257 2176}%
\special{pa 3263 2272}%
\special{pa 3263 2304}%
\special{pa 3264 2336}%
\special{pa 3260 2464}%
\special{pa 3258 2497}%
\special{pa 3256 2529}%
\special{pa 3253 2561}%
\special{pa 3251 2594}%
\special{pa 3245 2658}%
\special{pa 3242 2691}%
\special{pa 3240 2710}%
\special{fp}%
%
\special{pn 13}%
\special{pa 3120 760}%
\special{pa 3138 789}%
\special{pa 3157 819}%
\special{pa 3175 848}%
\special{pa 3192 877}%
\special{pa 3210 907}%
\special{pa 3244 965}%
\special{pa 3260 995}%
\special{pa 3275 1024}%
\special{pa 3290 1054}%
\special{pa 3304 1083}%
\special{pa 3317 1112}%
\special{pa 3329 1142}%
\special{pa 3340 1171}%
\special{pa 3350 1201}%
\special{pa 3358 1230}%
\special{pa 3365 1260}%
\special{pa 3371 1289}%
\special{pa 3376 1319}%
\special{pa 3379 1349}%
\special{pa 3380 1378}%
\special{pa 3380 1408}%
\special{pa 3378 1438}%
\special{pa 3375 1467}%
\special{pa 3370 1497}%
\special{pa 3364 1527}%
\special{pa 3356 1556}%
\special{pa 3338 1616}%
\special{pa 3327 1646}%
\special{pa 3315 1676}%
\special{pa 3302 1706}%
\special{pa 3289 1735}%
\special{pa 3244 1825}%
\special{pa 3227 1855}%
\special{pa 3211 1885}%
\special{pa 3194 1915}%
\special{pa 3176 1945}%
\special{pa 3159 1975}%
\special{pa 3120 2040}%
\special{fp}%
\put(20.1000,-19.3000){\makebox(0,0)[lb]{$\Theta_{1,1}$}}%
\put(19.9000,-22.5000){\makebox(0,0)[lb]{$\Theta_{1,0}$}}%
\put(26.6000,-24.9000){\makebox(0,0)[lb]{$\Theta_{2,0}$}}%
\put(26.4000,-19.9000){\makebox(0,0)[lb]{$\Theta_{2,1}$}}%
%
\special{pn 13}%
\special{pa 3710 750}%
\special{pa 3728 779}%
\special{pa 3747 809}%
\special{pa 3765 838}%
\special{pa 3782 867}%
\special{pa 3800 897}%
\special{pa 3834 955}%
\special{pa 3850 985}%
\special{pa 3865 1014}%
\special{pa 3880 1044}%
\special{pa 3894 1073}%
\special{pa 3907 1102}%
\special{pa 3919 1132}%
\special{pa 3930 1161}%
\special{pa 3940 1191}%
\special{pa 3948 1220}%
\special{pa 3955 1250}%
\special{pa 3961 1279}%
\special{pa 3966 1309}%
\special{pa 3969 1339}%
\special{pa 3970 1368}%
\special{pa 3970 1398}%
\special{pa 3968 1428}%
\special{pa 3965 1457}%
\special{pa 3960 1487}%
\special{pa 3954 1517}%
\special{pa 3946 1546}%
\special{pa 3928 1606}%
\special{pa 3917 1636}%
\special{pa 3905 1666}%
\special{pa 3892 1696}%
\special{pa 3879 1725}%
\special{pa 3834 1815}%
\special{pa 3817 1845}%
\special{pa 3801 1875}%
\special{pa 3784 1905}%
\special{pa 3766 1935}%
\special{pa 3749 1965}%
\special{pa 3710 2030}%
\special{fp}%
%
\special{pn 13}%
\special{pa 3880 750}%
\special{pa 3861 779}%
\special{pa 3841 809}%
\special{pa 3803 867}%
\special{pa 3785 896}%
\special{pa 3767 926}%
\special{pa 3749 955}%
\special{pa 3732 984}%
\special{pa 3716 1014}%
\special{pa 3701 1043}%
\special{pa 3687 1073}%
\special{pa 3674 1102}%
\special{pa 3662 1131}%
\special{pa 3651 1161}%
\special{pa 3642 1190}%
\special{pa 3634 1220}%
\special{pa 3628 1250}%
\special{pa 3623 1279}%
\special{pa 3621 1309}%
\special{pa 3620 1339}%
\special{pa 3620 1368}%
\special{pa 3623 1398}%
\special{pa 3628 1428}%
\special{pa 3634 1458}%
\special{pa 3642 1488}%
\special{pa 3651 1518}%
\special{pa 3661 1548}%
\special{pa 3672 1578}%
\special{pa 3684 1608}%
\special{pa 3696 1639}%
\special{pa 3722 1699}%
\special{pa 3748 1761}%
\special{pa 3760 1791}%
\special{pa 3773 1822}%
\special{pa 3795 1884}%
\special{pa 3805 1915}%
\special{pa 3813 1947}%
\special{pa 3821 1978}%
\special{pa 3828 2009}%
\special{pa 3834 2041}%
\special{pa 3839 2072}%
\special{pa 3847 2136}%
\special{pa 3853 2232}%
\special{pa 3853 2264}%
\special{pa 3854 2296}%
\special{pa 3850 2424}%
\special{pa 3848 2457}%
\special{pa 3846 2489}%
\special{pa 3843 2521}%
\special{pa 3841 2554}%
\special{pa 3835 2618}%
\special{pa 3832 2651}%
\special{pa 3830 2670}%
\special{fp}%
\put(39.4000,-19.6000){\makebox(0,0)[lb]{$\Theta_{4,1}$}}%
\put(39.0000,-23.5000){\makebox(0,0)[lb]{$\Theta_{4,0}$}}%
%
\special{pn 13}%
\special{pa 1260 1070}%
\special{pa 1780 1070}%
\special{fp}%
\put(46.3000,-10.3000){\makebox(0,0)[lb]{$s_{P_1}$}}%
%
\special{pn 13}%
\special{pa 1270 1595}%
\special{pa 1750 1595}%
\special{fp}%
%
\special{pn 13}%
\special{pa 1865 1590}%
\special{pa 2410 1590}%
\special{fp}%
%
\special{pn 13}%
\special{pa 2505 1590}%
\special{pa 3045 1595}%
\special{fp}%
%
\special{pn 13}%
\special{pa 3145 1590}%
\special{pa 3910 1590}%
\special{fp}%
%
\special{pn 13}%
\special{pa 3990 1595}%
\special{pa 4395 1595}%
\special{fp}%
%
\special{pn 13}%
\special{pa 4515 1600}%
\special{pa 4895 1595}%
\special{fp}%
\put(14.8500,-15.7000){\makebox(0,0)[lb]{$s_{P_{\tau}}$}}%
%
\special{pn 13}%
\special{pa 1250 1240}%
\special{pa 2000 1240}%
\special{fp}%
\special{pa 2160 1240}%
\special{pa 2360 1240}%
\special{fp}%
\special{pa 2520 1240}%
\special{pa 2970 1240}%
\special{fp}%
\special{pa 2970 1240}%
\special{pa 3280 1240}%
\special{fp}%
\special{pa 3440 1240}%
\special{pa 3540 1240}%
\special{fp}%
\special{pa 3710 1240}%
\special{pa 5170 1240}%
\special{fp}%
%
\special{pn 13}%
\special{pa 1230 1350}%
\special{pa 2000 1350}%
\special{fp}%
\special{pa 2160 1350}%
\special{pa 2690 1350}%
\special{fp}%
\special{pa 2810 1350}%
\special{pa 2950 1350}%
\special{fp}%
\special{pa 3090 1350}%
\special{pa 3530 1350}%
\special{fp}%
\special{pa 3710 1350}%
\special{pa 5180 1350}%
\special{fp}%
%
\special{pn 13}%
\special{pa 1900 1070}%
\special{pa 2600 1070}%
\special{fp}%
\special{pa 2790 1070}%
\special{pa 3210 1070}%
\special{fp}%
\special{pa 3400 1070}%
\special{pa 3600 1070}%
\special{fp}%
\special{pa 3780 1070}%
\special{pa 5180 1070}%
\special{fp}%
\end{picture}}%

Figure 2
\end{center}

\begin{lem}\label{lem:image-2}{
\begin{enumerate}
 \item[(i)]  The image of the fixed locus of $[-1]_{\varphi}$ is a nodal cubic $\mcE$ and a transversal line $\mcL_o$ to $\mcE$. Moreover, $\overline {q}\circ f(s_{P_{\tau}}) = \mcL_o$ and the node is the image of $\Theta_{4,1}$.
 
 \item[(ii)]  $\{f_{\mcQ, z_o}(\Theta_{1,1}), f_{\mcQ, z_o}(\Theta_{2,1}), f_{\mcQ, z_o}(\Theta_{3,1})\} =
 \mcE\cap \mcL_o$. We denote $p_i = f_{\mcQ, z_o}(\Theta_{i,1})$
 
\item[(iii)]  For $\{i,j,k\}=\{1,2,3\}$, $f_{\mcQ, z_o}(s_{P_{i,j}^\pm})$ ($i\not=j$) passes through $p_k$, $k\not=i,j$ ($i = 1, 2, 3$), respectively.

 \item[(iv)] Put $\mcL_k^\pm := f_{\mcQ, z_o}(s_{P_{i,j}^\pm})$.  Then   either
 (a) $\mcL_k^\pm$ is tangent to $\mcE$ at a point distinct to $p_k$  or 
 $p_k$  is an inflection point of $\mcE$ and $\mcL_k^\pm$ is
 an inflectional tangent line.
 
 \end{enumerate}
 }
 \end{lem}

\proof We prove the part about the node in statement (i). Since $s_{p_\tau}\Theta_{4,0}=1$, the remaining component(s) of the ramification locus must intersect $\Theta_{4,1}$ at two distinct points. Hence the image of $f(\Theta_{4,1})$ gives rise to a node on  the branch locus that is not on $\mcL_o$. The remaining statements can be proved in a similar way as that in Lemma \ref{lem:image-1}.  
\qed

Finally, as in the Introduction, we put 
\begin{eqnarray*}
P_4:=P_{1,3}^+=P_1\dot+P_3\dot+P_\tau \quad\quad
P_5:=P_{1,2}^+=P_1\dot+P_2\dot+P_\tau\\
P_6:=P_{2,3}^+=P_2\dot+P_3\dot+P_\tau \quad\quad
P_7:=P_{3,1}^-:=P_3\dot-P_1\dot+P_\tau
\end{eqnarray*}


Now let
\[
\mcB^1=\mcE +\mcL_0+\sum_{i=4}^6 \mcL_i, \quad
\mcB^2=\mcE +\mcL_0+\sum_{i=5}^7 \mcL_i.
\]
Then by  \cite[Theorem 3.2]{tokunaga14} and  Lemma~\ref{lem:key-comb1}, a $D_{2p}$-cover branched at $2(\mcE+\mcL_0)+p(\sum_{i=5}^7 \mcL_i)$ exists, but does not exist for $2(\mcE +\mcL_0)+p(\sum_{i=4}^6 \mcL_i)$. Hence, $(\mcB^1, \mcB^2)$  is a Zariski pair
if their combinatorics are the same.
 \bigskip 
%
%
%
%
\bigskip

{\bf Proof for Combinatorics 2}
We first note that our proof is similar to that of \cite[Theorem 5 (ii)]{tokunaga14}. Let $\varphi_{\mcQ, z_o} : S_{\mcQ, z_o} \to \PP^1$
 be the rational elliptic surface corresponding to Combinatorics 1-(b).  Under the same nation as before, Consider $[2]P_i$ ($i = 1, 2,3$)
 and their corresponding sections $s_{[2]P_i}$ $(i = 1, 2, 3)$, respectively. Since $\langle [2]P_i, [2]P_i \rangle = 2$ and
 $\gamma_{\mcQ, z_o}(s_{[2]P_i}) = (0, 0, 0, 0, 0)$,  we infer that $f_{\mcQ, z_o}(s_{[2]P_i})$ $(i = 1, 2, 3)$ are contact conics $\mcC_i$ ($i = 1, 2, 3)$  to
 $\mcQ$ and $z_o$ is one of the tangent points between $\mcC_i$ and $\mcQ$ by the argument in the proof of \cite[Theorem 5 (ii), p. 633]
 {tokunaga14}.  Now, for example, let  $\mcL = f_{\mcQ, z_o}(s_{P_1})$ and consider two curves
 \[
 \mcB^1 := \mcQ + \mcL + \mcC_1, \quad \mcB^2 := \mcQ + \mcL + \mcC_2.
 \]
 
 If there  exists a homeomorphism $h: (\PP^2, \mcB^1) \to (\PP^2, \mcB^2)$, it satisfies $h(\mcQ) = \mcQ$. Hence
 by \cite[Proposition~4.4]{tokunaga14}, there exist no homeomorphisms $(\PP^2, \mcB^1) \to (\PP^2, \mcB^2)$.
 Moreover, if both $\mcB^1$ and $\mcB^2$ have the Combinatorics 2, then $(\mcB^1, \mcB^2)$ is a Zariski pair. In \S 5, we show
 that such an example exists.


\section{Zariski triple and 4-ple for cubic-conic-line arrangements}\label{sec:3}

In \cite{bannai-tokunaga}, we give examples of Zariski $N$-ples for conic, 
conic-quartic arrangements. By similar  arguments to those in \cite{bannai-tokunaga},
we give a Zarisiki triple and $4$-ple for cubic-conic-line arrangements. 
Throughout this section, we use the terminology, notation  and results in \cite{bannai-tokunaga}, freely.
The combinatorics considered in this section is as follows:

\bigskip

{\bf Combinatorics 3.}
Let $\mcE$, $\mcL_o$ and $\mcC_i$ ($i = 1, 2, 3$) be as below. Put  $\mcB = \mcE  + \mcL_o + \sum_{i=1}^3 \mcC_i$:

\begin{enumerate}

\item[(i)] $\mcE$: (a)  a smooth cubic or  (b) a nodal cubic. 

\item[(ii)] $\mcL_o$: a transversal line to $\mcE$ and we put $\mcE \cap \mcL_o = \{p_1, p_2, p_3\}$.

\item[(iii)] $\mcC_i$ ($i = 1, 2, 3$):  contact conics to $\mcE + \mcL_o$. Each of them
is tangent to $\mcE + \mcL_o$ at four points

\item[(iv)] $\mcC_i$ and $\mcC_j$ intersect transversally for $i < j$ and  $\mcC_1\cap \mcC_2 \cap \mcC_3 = \emptyset$.

\end{enumerate}

%
%
%
%
%
%
%
%
%
%

%
%
%
%
%
%
%
%
%

 In the construction of plane curves with Combinatorics 3-(a) and (b),  how to find a contact conic $\mcC$ to 
 $\mcE + \mcL_o$ is crucial and we make use of bisections of elliptic surfaces as we did in \cite{bannai-tokunaga}.
 Let us recall that  a bisection is defined as follows:

 \begin{defin}\label{def:bisection}{\rm
 Let $\varphi : S \to C$ be an elliptic surface over $C$. Let $F$ be a general fiber of $\varphi$.
 A bisection of $\varphi$ is a horizontal curve $D$ with $FD = 2$. Here, a horizontal curve
 with respect to $\varphi$ is a curve that does not contain any fiber components.
 }
 \end{defin}

  Put $\mcQ = \mcE + \mcL_o$ and let $\varphi_{\mcQ, z_o} : S_{\mcQ, z_o} 
 \to \PP^1$ be the rational elliptic surface as before and let $f_{\mcQ, z_o} : S_{\mcQ, z_o} \to
 \PP^2$ be the generically $2$-$1$ morphism. Likewise \cite{bannai-tokunaga}, we construct
 a contact conic $\mcC$ as the image $f_{\mcQ, z_o}(D)$ of an irreducible bisection $D$.

{\sl Combinatorics 3-(a).} Our proof is almost parallel to that of \cite[Proposition 5]{bannai-tokunaga}.
Let $P_i$ ($i = 1, 2, 3, 4$) be the rational points introduced
in Proof for Combinatorics 1-(a) of  the previous section. Define $Q_i \in E_{\mcQ, z_o}(\CC(t))$ 
($i = 1, 2, 3, 4$)
by 
\[
[Q_0\dot{+}P_{\tau}\,\,  Q_1\,\,  Q_2\,\,  Q_3] = [P_0\,\,  P_1\,\,  P_2 \,\,  P_3] 
 \left [ \begin{array}{cccc}
         2 & -1 & -1 & -1 \\
        -1 & 2 & 0 & 0 \\
        -1 & 0 & 2 & 0 \\
        -1 & 0 & 0 & 2
          \end{array}
          \right ].
\]

Likewise \cite[p. 234]{bannai-tokunaga}, we now consider  six irreducible bisections  $D_0, \ldots, D_5$ as
follows:
\begin{enumerate}      
\item[(i)] $s(D_0) = -s_{Q_0}$,  $s(D_i) = - s_{Q_1} (1 \le i \le 3)$, $s(D_4) = - s_{Q_2}$, 
$s(D_5) = -s_{Q_0 + Q_1}$. Here $s(D_i)$ is the section determined uniquely by $D_i$ in \cite[Lemma~5.1]{shioda90}. 
\item[(ii)] $f_{\mcQ, z_o}(D_i)$ $(i = 0, 1,\ldots, 5)$  are contact conics to $\mcQ$ tangent at
$4$ distinct points.

\item[(iii)]  $\mcC_i$ and $C_j$ intersect transversally if $i \neq j$ and $\mcC_1\cap \mcC_2 \cap \mcC_3 = \emptyset$.
\end{enumerate}

Now if we put
\[
\mcB^1  :=  \mcQ + \mcC_1 + \mcC_2 + \mcC_3,  \,\mcB^2  :=   \mcQ + \mcC_0 + \mcC_1 + \mcC_2,\, 
\mcB^3  :=  \mcQ + \mcC_0 + \mcC_1 + \mcC_4, \, \mcB^4  :=   \mcQ + \mcC_0 + \mcC_1 + \mcC_5.
\]

Then by \cite[Theorem~4, Corollary~3, Corollary~4]{bannai-tokunaga}, we have

\begin{prop}\label{prop:key}{
\begin{enumerate}
\item[(i)] $\Cov_b(\PP^2, 2\mcQ + p(\mcC_i + \mcC_j), D_{2p}) \neq \emptyset$ if and only if $\{i, j\} \subset \{1, 2, 3\}$.
\item[(ii)]  $\Cov_b(\PP^2, 2\mcQ + p(\mcC_i + \mcC_j +\mcC_k), D_{2p}) \neq \emptyset$ if and only if $\{i, j, k\} \subset \{1, 2, 3\}$ or $\{0, 1, 5\}$.
\end{enumerate}
}
\end{prop}

From Proposition~\ref{prop:key},  we have

\begin{prop}\label{prop:z-4plet}{If $D_0, \ldots, D_5$ as above exist for 
$\mcQ$, $(\mcB^1, \mcB^2, \mcB^3, \mcB^4)$  is a Zarisiki 4-ple.
}
\end{prop}

%


{\sl  Combnatorics 3-(b).} Our proof is almost parallel to that of \cite[Thoerem 1]{bannai-tokunaga}.
Let $P_i$ ($i = 1, 2, 3$) be the rational points intruduced
in Proof for Combinatorics 1-(a) of  the previous section. Define $Q_i \in E_{\mcQ, z_o}(\CC(t))$ 
($i = 1, 2, 3$) by $Q_i := [2]P_i$ $(i = 1, 2, 3)$, respectively. We now consider $5$ bisections as
follows:

\begin{enumerate}      
\item[(i)]  $s(D_i) = - s_{Q_1} (1 \le i \le 3)$, $s(D_4) = - s_{Q_2}$, 
$s(D_5) = -s_{Q_3}$. Here $s(D_i)$ is the section determined uniquely by $D_i$ in \cite[Lemma~5.1]{shioda90}. 
\item[(ii)] $f_{\mcQ, z_o}(D_i)$ $(i =1,\ldots, 5)$  are contact conics to $\mcQ$ tangent at
$4$ distinct points.

\item[(iii)]  $\mcC_i$ and $C_j$ intersect transversally if $i \neq j$ and $\mcC_1\cap \mcC_2 \cap \mcC_3 = \emptyset$.
\end{enumerate}

Now if we put
\[
\begin{array}{ccccc}
\mcB^1  :=  \mcQ + \mcC_1 + \mcC_2 + \mcC_3,  & &
\mcB^2  :=   \mcQ + \mcC_1 + \mcC_2  + \mcC_4,  & &
\mcB^3 : =  \mcQ + \mcC_1 + \mcC_4 + \mcC_5. 
\end{array}
\]

Then by \cite[Theorem~4]{bannai-tokunaga}, we have

\begin{prop}\label{prop:key2}{
 $\Cov_b(\PP^2, 2\mcQ + p(\mcC_i + \mcC_j), D_{2p}) \neq \emptyset$ if and only if $\{i, j\} \subset \{1, 2, 3\}$.
}
\end{prop}

This shows that there exist no homeomrophism $(\PP^2, \mcB^i) \to (\PP^2, \mcB^j)$ if 
$i \neq j$ by a similar argument to the proof of \cite[Proposition~3]{bannai-tokunaga}. Thus we have

\begin{prop}\label{prop:z-triple}{If $D_1, \ldots, D_5$ as above exist for 
$\mcQ$, $(\mcB^1, \mcB^2, \mcB^3)$  is a Zarisiki triple.
}
\end{prop}


As we see in the next section, the six (resp. five) bisections as above exist for $\mcQ$ given
by an explicit equation. Thus we have

\begin{thm}\label{thm:z-tri-quar}{
 There exists a Zariski $4$-ple (resp. triple) for Comibnatorics 3-(a) (resp. (b)). 

%
%
%

}
\end{thm}

\begin{rem} {\rm (i) By increasing  the number of contact conics, our result can be generalized to
Zariski $N$-plet in the same way as in \cite{bannai-tokunaga}.

(ii) In Figure 2,  if we blow down $f(\Theta_{4,0})$, $f(O)$, $f(\Theta_{0,1})$, 
 $f(\Theta_{1,1})$,  $f(\Theta_{2,1})$ and $f(\Theta_{3,1})$, in this order, the the image 
 $\mcQ_1$of
 the branch locus $\Delta_f$ of $f$ consists of two smooth conics intersecting $4$ points. This 
 is the one we consider in \cite[Theorem1]{bannai-tokunaga}.

}
\end{rem}


\section{Examples}

We end this paper by giving explicit examples for Combinatorics $1$,  $2$ and $3$. 
Let us begin with explaining our method to construct explicit bisections briefly introduced in 
\cite{bannai-tokunaga}. We here use notation and terminology in \cite[2.2.3]{bannai-tokunaga} freely.

Let $U_i \cong \CC^2, \, (i = 1, 2)$ be affine open sets of $\Sigma_2$ with coordinates
$(t, x)$ $(i = 1)$ and $(s, x')$ $(i = 2)$ such that $t = 1/s, \, x = x'/s^2$. Suppose that $E_{\mcQ, z_o}$ is
given by a Weierstrass equation
\[
E_{\mcQ, z_o}: y^2  = f_{\mcT}(t, x), \quad f_{\mcT}(t, x) = x^3 + b_2(t) x^2 + b_3(t)x + b_4(t),
\]
where $b_i \in \CC[t]$, $\deg b_i \le i$ and $f_{\mcT}$ defines the trisection $\mcT$ on $U_i$. Let
$P = (x_P(t), y_P(t)) \in E_{\mcQ, z_o}(\CC(t))$. Consider the line in $\bbA^2_{\CC(t)}$ through
$P$ defined by
\[
L_P : y = l_P(t, x), \quad l_P(t, x) = r(t)(x - x_P(t)) + y_P(t), \, r(t) \in \CC(t).
\]
Then $f_{\mcT}(t, x) - l_P(t,x)^2$ factors into of the form
$(x - x_P(t))g(t, x)$, $g(t, x) \in \CC(t)[x]$. Suppose that $x_P(t), y_P(t) \in \CC[t]$ and we  choose
$r(t) \in \CC[t]$ such that $g(t, x) \in \CC[t,x]$ is irreducible and  the total degree of $g$ is $2$. Then the conic
$C(r(t), P)$ given by $g(t, x) = 0$ is a contact conic to $\mcQ$. Moreover, if we put
\[
f^*_{\mcQ, z_o}C(r(t), P) = C^+ + C^-,
\]
we have (i) $C^{\pm}$ are bisections and (ii) (if we choose $\pm$ suitably) $s(C^+)  = - s_P$.
We construct the bisections $D_i$ in the previous sections in this way. Now we go on to construct 
examples for each combinatorics.

\begin{exmple}[{\bf Combinatorics 1}]\label{ex:comb1} \rm

Let $[T, X, Z]$ be homogeneous coordinates of $\PP^2$ and let
 $(t, x):=(T/Z, X/Z)$ be affine coordinates for $\CC^{2} = \PP^{2} \setminus \{ Z = 0 \} $.

\medskip


\underline{{\sl Combinatorics 1-(a):}}  Consider $\mcE$ and $\mcL_o$ given by the affine equations:
\[
\mcE : x^{2} - t^{3} - 3{t}^{2} - 2t = 0, \,\,
{\mcL}_{o} : x = 0.
\]
$E_{\mcQ, z_o}$ is given by a Weierstarss equation
\[
E_{\mcQ, z_o} : y^2 = x(x^{2} - t^{3} - 3{t}^{2} - 2t)
\]

Put $p_{1}=[0,0,1]$, $p_{2}=[-1, 0, 1]$, and $p_{3}=[-2, 0, 1]$.  Choose $[0, 1, 0]$ as $z_{o}$. Then 
$\varphi_{\mcQ, z_o}$ has $4$ singular fibers of type $\III$ and $E_{\mcQ , z_{o}}(\CC (t)) \cong D_{4}^{\ast} \oplus \ZZ /2 \ZZ$. $P_{\tau}$ is given by $(0, 0)$. We choose  $P_{0}, P_{1}, P_{2}, P_{3}$ for a basis of $D_{4}^{\ast}$-part  as follows:
 \begin{itemize}
 \item $\displaystyle {P_{{0}}\, := \,\left [ \left( \frac{-1-i}4 \right)  \left( -t-1+i \right) ^{2},
 \frac {\sqrt{2-2\,i}}8\, \left( -t-1+i \right)  \left( t+1+i \right) ^{2} \right]}$,


\item  $\displaystyle{ P_{{1}}\, := \, \left [- \left(  \sqrt{2}-1 \right) t,t \left(  \sqrt{2}+t \right)  \sqrt{ \sqrt{2}-1}\right]}$.
\item $\displaystyle {P_{{2}}\, := \,\left [ \left( 1+i \right)  \left( t+1 \right) , \sqrt{-1-i} \left( t+1 \right) 
\mbox{} \left( -t-1+i \right) \right ]}$,

\item  $\displaystyle{ P_{{3}}\, := \, \left [-i \left(  \sqrt{2}+1 \right)  \left( t+2 \right) , \sqrt{i \left(  \sqrt{2}+1 \right) }
 \left( t+2+ \sqrt{2} \right)  \left( t+2 \right) \right ]}$,


\end{itemize}
 By straightforward 
 computation, we see that $\mcL_i: x - x_{P_i}(t) = 0$ are tangent lines for $\mcE$ through $p_i$ $(i = 1, 2, 3)$, while
 $x - x_{P_0}(t) = 0$ is a contact conic to $\mcQ$ through $z_o$,  and we infer that
 $\langle P_0, \, P_0 \rangle = 2$ and $\langle P_i , \, P_i \rangle = 1$.

 Also by explicit computation we have
 \[
 P_4: = P_2 \dot{-}P_3 \dot{+} P_{\tau} =
 \left [ -(\sqrt {2}+1)t, \, -\frac{\sqrt{-1-i}}{2} (i +1 +  {i}\sqrt{2})(-\sqrt{2}+t)t \right ].
 \]
  Let $\mcL_4: x - x_{P_4}(t) = 0$. Then since $\mcL_4$ is a tangent line, $\langle P_4, P_4\rangle=1$ from Remark \ref{rem:obs-1}, which implies $\langle P_2, P_3\rangle=\frac{1}{2}$. We can compute the height pairing for $ \langle P_i, P_j\rangle$ ($\{i,j\}\subset \{0,1,2,3\})$ in a similar way
 and we infer that $[\langle P_i, P_j \rangle ]$ is the $D_4^*$ given in \S 2.

 Then we see both $\mcQ + \sum_{i=1}^3\mcL_i$ and
$\mcQ + \sum_{i=2}^4$ have Combinatorics (1-a).

\medskip

 \begin{rem} {\rm The generators as above are computed by Ms. Emiko Yorisaki in her master's thesis 
 (\cite{yorisaki}).
} 
 \end{rem}


{\sl Combinatorics 1-(b)} Consider 
\begin{align*}
\mcE &: x^{2} - t^{3} - {t}^{2} = 0, \\
{\mcL}_{o} &: 2x - 3t - 3 = 0.
\end{align*}
Note that $\mcE$ has the node at $[0, 0, 1]$. Put $p_{1}=[-1, 0, 1]$, $p_{2}=[-3/4, 3/8, 1]$, and $p_{3}=[3, 6, 1]$, then $\mcE \cap {\mcL}_{o} = \{ p_{1}, p_{2}, p_{3} \}$. We put $\mcQ = \mcE + {\mcL}_{o}$ and choose $[0, 1, 0]$ as the distinguished point $z_{o}$. As we discussed in Section \ref{sec:3}, $E_{\mcQ , z_{o}}(\CC (t)) \cong (A_{1}^{\ast})^{\oplus 3}\oplus \ZZ /2 \ZZ$ and we can assume lines $\mcL_{i} \, (i=1, 2, 3)$ connecting the node of $\mcE$ and $p_{i} \, (i=1, 2, 3)$ are generators of $(A_{1}^{\ast})^{\oplus 3}$-part. Then we have the following coordinates:
\[
\begin{array}{cccc}
P_{\tau} = \left[ \frac{3t+3}{2}, 0 \right],  & 
 P_{1} = \left[ 0, \frac{\sqrt{6}t(t+1)}{2} \right],  &
 P_{2} =\left[ - \frac{t}{2}, \frac{\sqrt{2}t(4t+3)}{4} \right],  &
 P_{3} = \left[ 2t, \frac{\sqrt{-2}t(t-3)}{2} \right] 
 \end{array}
 \]

Furthermore, by computation on $E_{\mcQ, z_{o}}(\CC(t))$, we have $P_{4}$, $P_{5}$, $P_{6}$, and $P_{7}$ as follows:
\begin{align*}
P_4 & := P_3 \dot{+} P_1 \dot{+} P_{\tau} = \left[ {\frac { \left( 1-2 \sqrt {-3} \right)  \left( 26t+
21+3 \sqrt {-3} \right) }{52}}, \frac{\sqrt{2}(\sqrt{3}+\sqrt{-1})(4t+3)(2t+3+\sqrt{-3})}{16} \right] \\
P_5 & := P_1 \dot{+} P_2 \dot{+} P_{\tau} = \left[ -(1+\sqrt {3})( 2\,t-3-3\,\sqrt {3}),  \frac{\left( 2+\sqrt {3} \right)  \left( t-3 \right)  \left( t-6-4\,\sqrt {3} \right)
}{\sqrt{2}} \right] \\
P_6 & := P_2 \dot{+} P_3 \dot{+} P_{\tau} = \left[ -2 \sqrt{-1}(t+1), \frac{(2+ \sqrt{-1})(t+1)(t+2)}{\sqrt{2}} \right]\\
P_7 & := P_3 \dot{-} P_1 \dot{+} P_{\tau} = \left[ {\frac { \left( 1+2 \sqrt {-3} \right)  \left( 26t+
21-3 \sqrt {-3} \right) }{52}}, - \frac{\sqrt{2}(\sqrt{3}-\sqrt{-1})(4t+3)(2t+3-\sqrt{-3})}{16} \right]
\end{align*}
Let $f_{\mcQ, z_o}$ be as before and $\mcL_i := f_{\mcQ , z_{o}}(s_{P_{i}}) \, (i=4, 5, 6, 7)$.
Put
\[
\mcB^1 := \mcQ + \sum_{i=4}^6 \mcL_i,\,\,\, 
\mcB^2 := \mcQ + \sum_{i=5}^7 \mcL_i.
\]
Then we have a Zariski pair $(\mcB^1 , \mcB^2 )$.


\end{exmple}
\begin{exmple}[{\bf Combinatorics 2, 3-(b)}]\rm
We keep the same affine equations of $\mcE$, ${\mcL}_{o}$, ${\mcL}_{1}$, ${\mcL}_{2}$, and ${\mcL}_{3}$ as Example \ref{ex:comb1} 1-(b). Choose $[0, 1, 0]$ as the distinguished point $z_{o}$. As we discussed in Section \ref{sec:3}, we have a contact conic $\mcC$ to $\mcQ$ such that (i) $\mcC:= f_{\mcQ, z_o}(C^+)$ and (ii) $P_{C^+} = [2]P_1$. 

We put ${\mcB}^{i} := \mcQ + \mcC + {\mcL}_{i} \, (i = 1, 2, 3)$. From \ref{thm:comb-1}, we see that both of $({\mcB}^{1}, {\mcB}^{2})$ and $({\mcB}^{1}, {\mcB}^{3})$ are Zariski pairs having Combinatorics 2.

Next we will give explicit example of Zariski triple with Combinatorics 3-(b). Using the same basis and coordinates given  in Example \ref{ex:comb1} 1-(b), we obtain $\mcQ_1, \mcQ_2, \mcQ_3$ by explicit calculations as follows:
\begin{itemize}
\item $\mcQ_1=\displaystyle[2]P_{1}=\left[ \frac{{t}^{2}+9t+9}{6}, \frac{\sqrt {6}t \left( {t}^{2}-9\,t-9 \right)}{36}  \right]$
\item $\mcQ_2=\displaystyle[2]P_{2}=\left[\frac{1}{8}\left(t^2+16t+16 \right), \frac{\sqrt{2}}{32}\left(t+2\right)\left(t^2-16t-16\right) \right]$
\item $\mcQ_3=\displaystyle[2]P_{3}=\left[-\frac{1}{2}(t^2+t+1), \frac{\sqrt{2}\,i}{4} (t+2)(t^2-t-1)\right]$
\end{itemize}

We use the method given at the beginning of this section and construct bisections $D_1, \ldots, D_5$ and contact conics $\mcC_1, \ldots, \mcC_5$ corresponding to $-s_{\mcQ_1}$, $-s_{\mcQ_2}$, $-s_{\mcQ_3}$. 
The equations  can be calculated by using the data in the following table:
\begin{center}
\begin{tabular}{c|c|c}
conic & rational point & $r(t)$\\
\hline
$\mcC_j, (j=1,2,3)$ & $\mcQ_1$ & $\frac{1}{\sqrt{6}}t+b_j$\\
$\mcC_4$ & $\mcQ_2$ & $\frac{\sqrt{2}}{4}t+b_4$\\
$\mcC_5$ & $\mcQ_3$ & $-\frac{i}{\sqrt{2}}t+b_5$
\end{tabular}
\end{center}

The explicit equations for $C_{j} \, (j = 1, 2, 3), C_{4}$, and $C_{5}$ to $\mcE + {\mcL}_{o}$ become as follows:
\begin{align*}
\mcC _{j} : & {b_{j}}\,({b_{j}} + 6 \sqrt{6})\, t^2 -2 \sqrt{6}\, {b_{j}}\, tx + 6\, x^2 + 3\, {b_{j}}\, (3\,{b_{j}}+2 \sqrt{6})\, t - 6\, {b_{j}}^2 x + 9\, {b_{j}}^2=0, \\
\mcC _{4} : & (2\, {b_{4}}^2 +  30 \sqrt{2}\, {b_{4}} + 33)\, t^2 - 8 ( \sqrt{2}\, {b_{4}} -1)\, tx \\
&+ 16\, x^2 + 16 ( 2\, {b_{4}}^2 + 4 \sqrt{2}\, {b_{4}} + 3 )\, t - 8 ( 2\, {b_{4}}^2 -1)\, x + 16 ( 2\, {b_{4}}^2 + 2 \sqrt{2}\, {b_{4}} + 1 )=0, \\
\mcC _{5} : & ( {b_{5}}^2 - 6 )\, t^2 - 2 ( \sqrt{-2}\, {b_{5}} -2)\, tx - 2\, x^2 + ( {b_{5}}^2 - 4 \sqrt{-2}\, {b_{5}} - 6 )\, t \\
&+ 2 ( {b_{5}}^2 +2)\, x + {b_{5}}^2 - 2 \sqrt{-2}\, {b_{5}} -2=0.
\end{align*}
We put
\[
{\mcB}^{1}:= {\mcE} + {\mcL}_{o} + {\mcC}_{1} + {\mcC}_{2} + {\mcC}_{3}, \,
{\mcB}^{2}:= {\mcE} + {\mcL}_{o} + {\mcC}_{1} + {\mcC}_{2} + {\mcC}_{4}, \,
{\mcB}^{3}:= {\mcE} + {\mcL}_{o} + {\mcC}_{1} + {\mcC}_{4} + {\mcC}_{5}.
\]
If we choose general $b_{1}, b_{2}, b_{3}, b_{4}, b_{5} \in \CC$, then we have a Zariski triple $({\mcB}^{1}, {\mcB}^{2}, {\mcB}^{3})$.
\end{exmple}

\begin{exmple}[{\bf Combinatorics  3-(a)}] \label{ex:combi-3a}\rm

Let $P_0, P_1, P_2$ and $P_3 \in E_{\mcQ, z_o}(\CC(t))$ be the basis of $D^*_4$-part considered in Example~\ref{ex:comb1}, 1-(a).
Based on these points, we construct $Q_0, Q_1, Q_2, Q_3$ as in \S 3.  For $\mcQ_0, \mcQ_1, \mcQ_2$, we have explicit coordinates
\begin{align*}
&\bullet \displaystyle Q_{{0}}\, := \, \left[  - \left( \frac{\sqrt {2}+1}2 \right)  \left( 6-2\,t\sqrt {2}+{t}^{2}-4\,\sqrt {2}+4\,t \right)\right. , \\
 &\quad\quad\left.-\frac{1}{8}\,\sqrt {1-i} \left( -2+4\,i+3\,i\sqrt {2}-\sqrt {2} \right) 
 \left( 2\,t\sqrt {2}+{t}^{2}+4\,\sqrt {2}-2\,t-6 \right)  \left( -2+
\sqrt {2}-t \right) 
\right] \\
&\bullet \displaystyle Q_{{1}}\, := \, \left[  \frac{1}{4}(-1+i)\left( 2\,i+2\,it+2\,t+{t}^{2} \right) \right. ,\\
&\quad\quad\quad\left. - \frac{\sqrt{2}}{16}\left( i\sqrt {2-2\,i}+\,\sqrt {2-2\,i} \right)  \left( -i{t}^{2}+{t}^{3}-2\,it+3\,{t}^{2}+2-2\,i+4\,t \right) \right] \\
&\bullet \displaystyle Q_{{2}}\,:=\left[ \frac{1}{4}\left( 1+i\right)  \left( -t-1+i \right) ^{2}\right.,\left. \frac{1}{8}\left( -1+i\right) \sqrt {-1-i} \left( -t-1+i \right)  \left( t+1+i \right) ^{2} \right]
\end{align*}

We use the method given at the beginning of this section and construct bisections $D_0, \ldots, D_5$ and contact conics $\mcC_0, \ldots, \mcC_5$ corresponding to $-s_{\mcQ_0}$, $-s_{\mcQ_1}$, $-s_{\mcQ_2}$  and $-s_{\mcQ_0\dot+\mcQ_1}$. 
The equations  can be calculated by using the data in the following table. 
\begin{center}
\begin{tabular}{c|c|c}
conic & rational point & $r(t)$\\
\hline
$\mcC_0$ & $\mcQ_0$ & $\frac{1}{4}\sqrt{1-i}\left(-2i-i\sqrt{2}+\sqrt{2}\right) t+b_0$\\
$\mcC_j, (j=1,2,3)$ & $\mcQ_1$ &  $\frac{i}{2}\sqrt {1-i}\,t+b_j$\\
$\mcC_4$ & $\mcQ_2$ & $\-\frac{1}{4}\sqrt{2-2i}(i+1)t+b_4$ \\
$\mcC_5$ & $\mcQ_1\dot+\mcQ_2$ & $\frac{1}{4}(1-i)\left(i\sqrt {2}+i+\sqrt {2}\right) t+b_5$
\end{tabular}
\end{center}

For $\mcC_0$, by using $P=\mcQ_0$ and $r(t)=\left(-2i-i\sqrt{2}+\sqrt{2}\right)\frac{\sqrt{1-i}}{4} t+b_0$ as above, we have the explicit equation

\begin{align*}
&\mcC_0: \left( {t}^{2}+4\,t-2\,t\sqrt {2}+6-4\,\sqrt {2}+2\,\sqrt {2}x-2\,x
 \right) {b}^{2}\\
 &+\sqrt {1-i} \left( i\sqrt {2}{t}^{2}-i\sqrt {2}tx+4\,
i\sqrt {2}t-i{t}^{2}+2\,\sqrt {2}{t}^{2}-\sqrt {2}tx+4\,i\sqrt {2}-6\,
it+10\,t\sqrt {2}-3\,{t}^{2} \right.\\
&\left.\quad\quad+2\,tx-14-6\,i+10\,\sqrt {2}-14\,t
 \right)b \\
 &-\left( 7-5\,\sqrt {2} \right)  \left( 6\,\sqrt {2}tx-4\,{x}^{2}
\sqrt {2}+2\,t\sqrt {2}+2\,\sqrt {2}x-{t}^{2}+8\,tx-6\,{x}^{2}+2\,x-2
 \right) 
=0
\end{align*}
We omit the equations of the other conics as they are rather long. 
We put 
\[
\mcB^1  :=  \mcQ + \mcC_1 + \mcC_2 + \mcC_3,  \,\mcB^2  :=   \mcQ + \mcC_0 + \mcC_1 + \mcC_2,\, 
\mcB^3  :=  \mcQ + \mcC_0 + \mcC_1 + \mcC_4, \, \mcB^4  :=   \mcQ + \mcC_0 + \mcC_1 + \mcC_5.
\]
It can be checked that for a general choice of $b_0,\ldots, b_5$, 
these curves  have Combinatorics 3-(a), and we have a Zariski 4-ple.

\end{exmple}


%

\noindent Shinzo BANNAI\\
National Institute of Technology, Ibaraki College\\
866 Nakane, Hitachinaka-shi, Ibaraki-Ken 312-8508 JAPAN \\
{\tt sbannai@ge.ibaraki-ct.ac.jp}\\

\noindent Hiro-o TOKUNAGA, Momoko YAMAMOTO\\
Department of Mathematics and Information Sciences\\
Tokyo Metropolitan University\\
1-1 Minami-Ohsawa, Hachiohji 192-0397 JAPAN \\
{\tt tokunaga@tmu.ac.jp, yamamoto-momoko@ed.tmu.ac.jp}
%
%

{\tt }

\vspace{0.5cm}
      
 \end{document}